\documentclass[a4paper,oneside,11pt,reqno]{amsart}
\usepackage[ansinew]{inputenc}
\usepackage{amsfonts}
\usepackage{amssymb}
\usepackage{amsmath,amsthm,amsfonts,amssymb}
\usepackage[left=2.5cm,right=2.5cm,top=3cm,bottom=3cm]{geometry}
\usepackage[dvips]{graphicx,color}

\theoremstyle{plain}
\newtheorem{lem}{Lemma}[section]
\newtheorem{condition}{Condition}
\newtheorem{thm}[lem]{Theorem}
\newtheorem{prop}[lem]{Proposition}
\theoremstyle{definition}
\newtheorem{Def}[lem]{Definition}
\newtheorem{Rem}[lem]{Remark}

\newcommand{\ve}{V}

\DeclareMathOperator{\ess}{ess}

\DeclareMathOperator{\Span}{Span}

\newcommand{\bu}{\mathbf{u}}
\newcommand{\CZ}{\mathcal{Z}}
\newcommand{\bw}{\mathbf{w}}
\newcommand{\wie}{\widetilde{\eta}}

\newcommand{\bb}{\mathbf{B}}
\newcommand{\xm}{X^m}
\newcommand{\xms}{\xm(s)}
\newcommand{\prm}{\Pi_m}
\newcommand{\bbu}{\tilde{\bu}^m}
\newcommand{\bbus}{\bbu(s)}

\newcommand{\bv}{\mathbf{v}}

\newcommand{\bum}{\mathbf{u}^m}
\newcommand{\bums}{\bum(s)}
\newcommand{\buns}{\bum(s-)}

\newcommand{\bae}{{\mathbb{E}}}

\newcommand{\DEQS}{\begin{eqnarray*}}
\newcommand{\EEQS}{\end{eqnarray*}}
\newcommand{\DEQSZ}{\begin{eqnarray}}
\newcommand{\EEQSZ}{\end{eqnarray}}
\newcommand{\EE}{\mathbb{E}}

\newcommand{\FF}{\mathbb{F}}

\newcommand{\Rb}[1]{{\mathbb{R}_{#1}}}
\newcommand{\bcase}{\begin{cases}}
\newcommand{\ecase}{\end{cases}}

\newcommand\del[1]{}

\newcommand\cadlag{c{\`a}dl{\`a}g }

\newcommand\dela[1]{}%{\Blue{#1}}

\newcounter{gg11}
\newenvironment{list-a}
{\begin{list} {{\rm (\alph{gg11})}}
{\usecounter{gg11}
\setlength{\leftmargin}{0.3cm}
\setlength{\topsep}{0.1cm}
\setlength{\itemsep}{0.0cm}
\setlength{\parsep}{0.1cm}
\setlength{\itemindent}{0.0cm}
\setlength{\parskip}{0.0cm}}}
{\end{list}}
\newcounter{gg111}
\newenvironment{list-n}
{\begin{list} {{\rm (\roman{gg111})}}
{\usecounter{gg111}
\setlength{\leftmargin}{0.5cm}
\setlength{\topsep}{0.1cm}
\setlength{\itemsep}{0.0cm}
\setlength{\parsep}{0.1cm}
\setlength{\itemindent}{0.0cm}
\setlength{\parskip}{0.0cm}}}
{\end{list}}
\newcounter{gg1}

\newcounter{lil}

\def\eps{\varepsilon}
\newcommand{\lk}{\left}

\newcommand{\rk}{\right}
\newcommand{\la}{{\langle}}
\newcommand{\ra}{{\rangle}}

\newcommand{\LL}{{\rm I \kern -0.2em L}}

\newcommand{\be} {\begin{enumerate} }
\newcommand{\ee} {\end{enumerate} }

\newcommand{\CS}{{{ \mathcal S }}}

\newcommand{\CM}{{{ \mathcal M }}}
\newcommand{\CP}{{{ \mathcal P }}}

\newcommand{\RR}{{\mathbb{R}}}

\newcommand{\DD}{{\rm I \kern -0.2em D}}
\newcommand{\dd}{{\rm I \kern -0.2em D}}

\newcommand{\NN}{\mathbb{N}}

\newcommand{\PP}{{\mathbb{P}}}

\newcommand{\TT}{{\rm I \kern -0.2em T}}

\newcommand{\DEQ}{\begin{eqnarray}}
\newcommand{\EEQ}{\end{eqnarray}}

\begin{document}
\title[Stochastic evolution equations for Non-Newtonian fluids]
{On stochastic evolution equations for nonlinear bipolar fluids: well-posedness and some properties of the solution}
\author[E. Hausenblas]{Erika Hausenblas}
\email[E. Hausenblas]{erika.hausenblas@unileoben.ac.at}
\author[P. A. Razafimandimby]{Paul Andr\'e Razafimandimby}
\email[P.A. Razafimandimby]{paulrazafi@gmail.com}

\address{Department of Mathematics and Information Technology\\
Montan University Leoben\\ Franz Josef Str. 18, 8700 Leoben,
     Austria}
%\address[P. A. Razafimandimby]{Mathematics Section\\ The Abdus Salam International Center for Theoretical Physics\\
%Strada Costiera, 11, I - 34151 Trieste, Italy}
%\curraddr[P. A. Razafimandimby]{Department of Mathematics and
%Applied Mathematics\\
%University of Pretoria\\
%Lynwood Road, Pretoria 0002, South Africa}
%\address[M. Sango]{Department of Mathematics and
%Applied Mathematics\\
%University of Pretoria\\
%Lynwood Road, Pretoria 0002, South Africa}
\date{\today}

\begin{abstract}
We investigate the stochastic
evolution equations describing the motion of a  Non-Newtonian fluids excited
by multiplicative noise of L\'evy type. We show that the system we consider has a unique global strong solution.
We also give some results concerning the properties of the solution. Mainly we prove that the unique solution satisfies
the Markov-Feller property. This enables us to prove by means of some results form ergodic theory that the semigroup
 associated to the unique solution admits at least an invariant measure which is ergodic and tight on a subspace of the Lebesgue space
 $L^2(\mathcal{O})$.
\end{abstract}

\maketitle

\section{Introduction}

Turbulence in Hydrodynamics is one of the most fascinating and difficult problems in Mathematics and
in applied sciences in general. Many scientists believe that Newtonian law or the Navier-Stokes Equations (NSE for short) can
accurately decscribe the most intricate complexities of turbulence in fluids flows.
However, there are mainly two major obstacles for the mathematical study of turbulent flows. First, it is well known that
the question of whether the three dimensional NSE admits or not a
unique weak solution for all time still remains open. As it is not alwaays easy to prove
the existence of a global attractor in the case of lack of uniqueness of solution, this becomes a daunting obstacle for
the investigation of the long-time behaviour of the Navier-Stokes equations which is very important for a better
understanding of turbulence and some physical features of the fluids. We refer, for instance, to
\cite{BABIN+VISHIK}, \cite{Constantin+Foias+Temam}, \cite{Constantin+Foias+Temam-2}, \cite{Robinson}, and \cite{Temam-Infinite} for
some results in this direction. Second, there are a lot of fluid models exhibiting turbulent behavior that cannot be described
by the Navier-Stokes equations. To overcome these problems one generally has to use other model of fluids
 or some regularizations, which might be of mathematical nature, of the Newtonian law. This has motivated
  many scientists to consider fluids such that their stress tensor is a nonlinear functions of the strain rate.
  This class of fluids forms the family of Non-Newtonian fluids. One examples of such fluids is the nonlinear Bipolar fluids which are themselves contained
   in the class of multipolar fluids. The theory of viscous multipolar fluids was initiated by Necas and
Silhavy \cite{NECAS2}, and developed later on in numerous work of prominent scientists such as Necas, Novotny and Silhavy
\cite{NECAS1}, Bellout, Bloom and Necas \cite{BELLOUT2}. Although Bipolar fluids ressembles to the models that Ladyzhenskaya considered
in \cite{LADY1} and \cite{LADY2} they differ in two aspects. First both bipolar fluids and Ladyzhenskaya allow for a nonlinear velocity
 dependent viscosity, but in contrast to Bipolar fluids the Ladyzhenskaya models do not incorporate a higher-order velocity gradients.
 Second, in contrary to the Ladyzhenskaya models the theory of multipolar fluids is compatible with the basic principles
of thermodynamics such as the Clausius-Duhem inequality and the principle of
frame indifference. Moreover, results up to date indicate that the theory
of multipolar fluids may lead to a better understanding of hydrodynamic
turbulence (see for example \cite{BELLOUT3}).

Around the 70s Bensoussan and Temam \cite{Bensoussan-Temam} started the investigation of stochastic version of dynamical equations
for Newtonian turbulent fluids.
The Stochastic Partial Differential Equations they analyzed are obtained by
adding noise terms to the deterministic NSE. This approach is basically
motivated by Reynolds' work which stipulates that the velocity of a fluid
particle in turbulent regime is composed of slow (deterministic) and fast
(stochastic) components. While this belief was based on empirical and
experimental data, Rozovskii and Mikulevicius were able to derive the models
rigorously in their recent work \cite{ROZOVSKII1}, thereby confirming the
importance of this approach in hydrodynamic turbulence. It is also pointed out in some recent articles like
\cite{Flandoli2}  and \cite{KAP}) that some rigorous information on questions in
Turbulence might be obtained from stochastic versions of the equations of
fluid dynamics. Since the pioneering work of Bensoussan and Temam \cite{Bensoussan-Temam} on
stochastic Navier-Stokes equations, stochastic models for Newtonian fluid dynamics and SPDEs in general have been the
object of intense investigations which have generated several important
results. We refer, for instance, to \cite{Alb+Brz+Wu_2010}, \cite{bensoussan},%
\cite{Breckner}, \cite{Brzezniak3}, \cite{Brzezniak4}, \cite{ZB+EH+JZ}, \cite{Caraballo1},
\cite{Caraballo}, \cite{Caraballo3},  \cite{daprato},\cite{DAPRATO}, \cite{DEUGOUE2}, \cite{DongNS},\cite{ROZOVSKII1}, \cite{pardoux}, \cite{Peszat}%
,\cite{sango}, \cite{sango2}, \cite{Taniguchi}.  However, there are only very
few results for the dynamical behaviour of stochastic models for
Non-Newtonian fluids (see \cite{GUO1}, \cite{PAUL1}, \cite%
{PAUL2}, \cite{PAUL3}, \cite{EH+PAR}).

In this paper, we are
 interested in the L\'evy driven SPDEs for the nonlinear Bipolar fluids.
%whose reduced stress tensor $\hat{%
% \mathbb{T}}(\mathcal{E}(\bu))=\mathbf{T}(\mathcal{E}(\bu))$ is given by
% \begin{equation*}
% \mathbf{T}(\mathcal{E}(\bu))=2\kappa_0(1+ |\mathcal{E}(\bu)|^2)^{\frac{p-2}{2}}%
% \mathcal{E}(\bu)-2\kappa_1\Delta\mathcal{E}(\bu).
% \end{equation*}
More precisely, let $d=2,3$, and $\mathcal{O}\subset \mathbb{R}^d$ be a smooth bounded
open domain,  we consider
\begin{equation}  \label{1}
\begin{cases}
d\bu+\left[\bu\cdot\nabla \bu-\nabla\cdot\mathbf{T}(\mathcal{E}(\bu))+\nabla \pi\right]%
dt=\int_Z\sigma(t,\bu,z)\wie(dz,dt),\,\, x\in \mathcal{O}, t\in (0,T], \\
\bu(x,0)=\bu_0,\,\, x\in \mathcal{O}, \\
\nabla\cdot \bu=0,\,\, x\in \mathcal{O}, t\in [0,T], \\
\bu(x,t)=\tau_{ijl}n_jn_l=0,\,\, x\in \partial \mathcal{O}, t\in (0,T)%
\end{cases}%
\end{equation}
where $\bu$ is the velocity of the fluids, $\pi$ its pressure, $\mathbf{n}$ denotes the
normal exterior to the boundary and
\begin{equation*}
\mathcal{E}(\bu)=\frac{1}{2}\left(\frac{\partial \bu_i}{\partial x_j}+\frac{%
\partial \bu_j}{\partial x_i}\right), |\mathcal{E}(\bu)|^2=\sum_{i,j=1}^n|%
\mathcal{E}_{ij}(\bu)|^2,
\end{equation*}
\begin{equation*}
\mathbf{T}(\mathcal{E}(\bu))=2\kappa_0(\varkappa+|\mathcal{E}(\bu)|^2)^{\frac{p-2}{2}}%
\mathcal{E}(\bu)-2\kappa_1\Delta\mathcal{E}(\bu).
\end{equation*}
The quantities $\kappa_0$, $\kappa_1$ and $\varkappa$ denote positive constants.
Here $\wie$ is a compensated Poisson random measure  defined on a prescribed
probability space $(\Omega, \mathcal{F}, \mathbb{P})$ and taking its values
in a separable Hilbert space $H$ to be defined later. The system \eqref{1}
describes the equations of motion of isothermal incompressible nonlinear
bipolar fluids excited by random forces.

For $p=2,\,\, \kappa_1=0,\,\, \sigma\equiv 0$, \eqref{1} is the Navier-Stokes
equations which has been extensively studied (see, for instance, \cite{temam}%
). If $1<p<2$ then the fluid is shear thinning, and it is shear thickening
when $2<p$. The problem \eqref{1} is as interesting as the Navier-Stokes
equations. It contains two nonlinear terms which makes the problem as
difficult as any nonlinear evolution equations. During the last two decades,
the deterministic version of \eqref{1} has been the object of intense
mathematical investigation which has generated several important results. We
refer to \cite{BELLOUT1}, \cite{BELLOUT4}, \cite{BELLOUT5}, \cite{MALEK},
\cite{ROKYTA} for relevant examples.  Despite these numerous  results there
are still a lot of open problems related to the mathematical theory of
multipolar fluids. Some examples are the existence of weak solution for all
values of $p$, the uniqueness of such weak solutions and many more. We
refer, for instance, to \cite{BELLOUT4}, \cite{RUZICKA} and \cite{ROKYTA}
for some discussions about these challenges.

 For  $1<p$ and the noise is replaced by a cylindrical Wiener process, the existence of
martingale and stationary solution of \eqref{1} was established in \cite{GUO1}.
Razafimandimby and Sango studied the exponential stability and some stabilization of \eqref{1} with $1<p\le 2$ and with a Wiener noise.
It seems that this article is the first work studying the L\'evy driven SPDEs  \eqref{1}. Our first main goal is to prove
the existence and uniqueness of strong solution which should be understood in the sense of stochastic differential equations.
To achieve this goal we mainly follow the idea initially developed by Breckner in \cite{Breckner} (see also \cite{breckner}) and used later in many
articles such as  \cite{Caraballo}, \cite{DEUGOUE2}, \cite{PAUL-13}. This method is based on Galerkin approximation and
it allows to prove that the whole sequence of the Galerkin approximation converges in mean square to the exact solution.
The second of the present paper is to give some partial
results concerning the properties of the solution. We concentrate on proving that the solution satisfies the Markov-Feller property
 which enables us to prove that \eqref{1} admits at least an invariant measure which is ergodic and tight on a subspace of the Lebesgue
 space $L^2(\mathcal{O})$. Unfortunately we could not proceed further and prove the uniqueness of the invariant measure. The investgation of the uniqueness of ergodic SPDEs driven by
 pure jump noise seems to be very difficult and out of reach of the most recent methods used to prove the uniqueness of invariant measure of SPDEs.
 We postpone this investigation in future work.

 To close this introduction we give the outline of the article. In Section 2 we give most of the notations and necessary
 preliminary used throughout the work. By means of Galerkin approximation we show the existence of strong solution in Section 3.
 The pathwise uniqueness of the solution and the convergence of the whole sequence of Galerkin approximate solution to the exact solution are proved
 in Section 4. Section 5 is devoted to the investigation of some properties of the strong solution.

\noindent \textbf{Notations}:

\noindent \textit{By $\mathbb{N}$ we denote the set of nonnegative
integers, i.e. $\mathbb{N}=\{0,1,2,\cdots\}$ and by
$\bar{\mathbb{N}}$ we denote the set $\mathbb{N}\cup\{+\infty\}$.
Whenever we speak about $\mathbb{N}$ (or
$\bar{\mathbb{N}}$)-valued measurable functions we implicitly
assume that the  set is equipped with the trivial $\sigma$-field
$2^\mathbb{N}$ (or $2^{\bar{\mathbb{N}}}$). By $\Rb{+}$ we will
denote the  interval $[0,\infty)$ and by $\Rb{\ast}$ the set
$\mathbb{R}\setminus\{0\}$. If $X$ is a topological space, then by
$\mathcal{B}(X)$ we will denote the Borel $\sigma$-field on $X$.
By $\lambda_d$ we will denote the Lebesgue measure on
$(\mathbb{R}^d,\mathcal{B}(\mathbb{R}^d))$, by $\lambda$  the
Lebesgue measure on $(\mathbb{R},\mathcal{B}(\mathbb{R}))$.}

\textit{If $(S,\CS)$ is a measurable space then by  $M(S)$ we
denote the set of all real valued  measures on $(S,\CS)$, and  by
$\CM(S)$ the $\sigma$-field
 on $M(S)$ generated by functions
$i_B:M(S) \ni\mu \mapsto \mu(B)\in \RR$, $B\in \CS$. By $M_+(S)$
we denote the set of all non negative measures on $S$,  and  by
$\CM(S)$ the $\sigma$-field
 on $M_+(S)$ generated by functions
$i_B:M_+(S) \ni\mu \mapsto \mu(B)\in \RR_+$, $B\in \CS$. Finally,
by $M_I( S)$ we denote the family of all
$\overline{\mathbb{N}}$-valued measures on $(S,\CS)$,   and  by
$\CM_I(S)$ the $\sigma$-field
 on $M_I(S)$ generated by functions
$i_B:M(S) \ni\mu \mapsto \mu(B)\in \bar\NN$, $B\in \CS$. \del{ and
by $\CM_I(S\times \RR_+)$  the $\sigma$-field on $M_I(S\times
\RR_+)$ generated by functions $i_B:M_I(S\times \RR_+) \ni\mu
\mapsto \mu(B)\in \overline{\mathbb{N}}$, $B\in \CS$.}}
%
%
%If $(S,\CS)$ is a measurable space, then by  $\CM^+(S)$ we denote the set of all positive measures on $S$,
\textit{If $(S,\CS)$ is a measurable space then  we will denote
by $\CS\otimes \mathcal{B}({\mathbb{R}_+})$  the product
$\sigma$-field on $S\times \mathbb{R}_+$ and by $\nu\otimes
\lambda$ the product measure of $\nu$ and the Lebesgue measure
$\lambda$.}

\section{Mathematical settings of the problem \eqref{1}}

Throughout this paper we mainly use the same notations as in \cite{GUO1}. By
$L^q(\mathcal{O})$ we denote the Lebesgue space of $q$-th integrable functions with
norm $||\cdot||_{L^q}$. For the particular case $q=2$, we denote its norm by
$||\cdot||$. For $q=\infty$ the norm is defined by $||\bu||_{L^\infty}=%
\ess \sup_{x\in D}|\bu(x)|$, where $|x|$ is the {{Euclidean}} norm of the vector $%
x\in \mathbb{R}^n$. The Sobolev space $\lbrace \bu\in L^q(\mathcal{O}): D^k \bu\in L^q(\mathcal{O}),
k\le \sigma \rbrace$ with norm $||\cdot||_{q,\sigma}$ is denoted by $%
W^{q,\sigma}(\mathcal{O})$. $C(I,X)$ is the space of continuous functions from the
interval $I=[0,T]$ to $X$, and $L^q(I,X)$ is the space of all measurable
functions $u:[0,T]\rightarrow X$, with the norm defined by
\begin{equation*}
||\bu||^q_{L^q(I,X)}=\int_0^T||\bu(t)||_X^q dt, \,\, q\in [1,\infty),
\end{equation*}
and when $q=\infty$, $||\bu||_{L^\infty(I,X)}=\ess\sup_{t\in [0,T]}||\bu(t)||_X$.

The
mathematical expectation associated to the probability space $(\Omega,
\mathcal{F}, \mathbb{P})$  is denoted by $\mathbb{E}$ and as above we also define the space $L^q(\Omega, X)$.

We proceed with the definitions of some additional spaces frequently used in
this work. We define a space of smooth functions with support
strictly contained in $\mathcal{O}$ and satisfying the divergence free condition:
\begin{equation*}
\mathcal{V}=\lbrace \bu \in \mathcal{C}_{c}^{\infty }(\mathcal{O}):\nabla\cdot
\bu=0\rbrace.
\end{equation*}
We denote by $H$ the closure of $\mathcal{V}$ with norm $|\cdot|$ in $%
L^2(\mathcal{O})$. It is a Hilbert space when equipped with the $L^2(\mathcal{O})-$inner product
$(.,.)$. $\mathcal{H}^\sigma$ is the closure of $\mathcal{V}$ in $W^{2,\sigma}(\mathcal{O})$
with the norm $||\cdot||_{2,\sigma}$. We denote by $\lVert \bu\rVert_\sigma$ the norm induced by $\lVert \bu\rVert_{2,\sigma}$ on $\mathcal{H}^\sigma$.
We denote by $\mathcal{H}^{-\sigma}$ the dual space  of $\mathcal{H}^\sigma$ ($%
\sigma\ge 1$) wrt the norm $\lVert\bu\rVert_\sigma$. If $\sigma=2$, then $V=\mathcal{H}^2$ and $V^\ast$ is the dual space
of $V$. The duality product between $V$ and $V^\ast$ is denoted by $%
\langle.,.\rangle$. {{It should be noted that $V$ is not the usual
space of divergence-free functions of $W^{2,1}(\mathcal{O})$ used for the Navier-Stokes
equations. Here it is a space of divergence-free functions of $W^{2,2}(\mathcal{O})$}}. We
assume throughout the paper that there exists a positive constant $\lambda_1$ such
that that the Poincar\'e inequalities type
\begin{equation}  \label{POINCARE}
||\bu||^2_\sigma \le \frac{1}{\lambda_1}||\bu||^2_{\sigma +1}, \forall \bu\in \mathcal{H}^{\sigma +1}, \sigma\ge 0,
\end{equation}
hold.

% Note that $H$ is a separable Hilbert space,  so we can associate with it an
% orthonormal basis $\lbrace e_i:i=1,2,...\rbrace$.

As mentioned in the introduction we will study a stochastic model for
a nonlinear bipolar fluids excited by random forces. In the following
lines we describe the forces acting on the fluids. Let $(Z,\mathcal{Z})$ be a separable metric space and let
 $\nu$ be a $\sigma$-finite positive measure on it. Suppose that $\mathfrak{P}=(\Omega,\mathcal{F},\FF,\mathbb{P})$ is a filtered
probability space, where $\FF=(\mathcal{F}_t)_{t\geq 0}$ is a filtration,
and $\eta: \Omega\times \mathcal{B}(\mathbb{R}_+)\times\mathcal{Z}\to
\bar{\mathbb{N}}$ is a time homogeneous Poisson random measure with
the intensity measure $\nu$ defined over the filtered probability
space $\mathfrak{P}$. We will denote by $\tilde\eta=\eta-\gamma$
the compensated Poisson random measure associated to $\eta$ where
the compensator $\gamma$ is given by
$$
\mathcal{B}(\mathbb{R}_+)\times \mathcal{Z}\ni (A,I)\mapsto  \gamma(A,I)=\nu(A)\lambda(I)\in\mathbb{R} _ +.
$$

We assume that $p\in (1,2]$. This will be fixed in the whole section.
% It is then known, see e.g. \cite{Brz+Haus_2009}, \cite{Zhu_2010}, that there exists a unique
% continuous linear operator $I$ which associates to each progressively measurable
% process $\xi:\Rb{+}\times Z\times \Omega \to \ve$ such that
{Let
$\mathcal{M}^2(\Rb{+},L^2(Z,\nu,H))$ be the class of all progressively measurable processes
$\xi:\Rb{+}\times Z\times \Omega \to \ve$ satisfying the condition
\begin{equation}\label{cond-2.01}
          \mathbb{E} \int_0^T\int_Z\vert \xi(r,z)\vert^2_H\nu(dz)\,dr <\infty, \quad T>0.
\end{equation}
If $T>0$, the class of all
progressively measurable processes $\xi:[0,T]\times  Z\times
\Omega \to \ve$ satisfying the condition \eqref{cond-2.01} just
for this one $T$,  will be denoted by
$\mathcal{M}^2(0,T,L^2(Z,\nu,H))$.
Also, let $\mathcal{M}_{step}^2(\Rb{+},L^2(Z,\nu,H))$ be the space of all processes $\xi \in \mathcal{M}^2(\Rb{+},L^2(Z,\nu,H))$ such that
$$
\xi(r) = \sum_{j=1} ^n 1_{(t_{j-1}, t_{j}]}(r) \xi_j,\quad 0\le r,
$$
where $\{0=t_0<t_1<\ldots<t_n<\infty\}$ is a partition of $[0,\infty)$, and
for all $j$,    $\xi_j$ is  an $\mathcal{F}_{t_{j-1}}$ measurable random variable. For any $\xi\in\mathcal{M}_{step}^2(\Rb{+},L^2(Z,\nu,H))$ we set
\begin{equation}
\label{eqn-2.02} \tilde{I}(\xi)= % \int_0^t \int_Z \xi(r,z) \tilde{\eta}(dz,dr)
\sum_{j=1}^n  \int_Z  \xi_j (z) \tilde\eta \left(dz,(t_{j-1}, t_{j}] \right).
\end{equation}
Basically, this is the definition of stochastic integral of a random step process $\xi$ with respect to the compound random Poisson measure $\tilde{\eta}$.
The extension of this integral on $\mathcal{M}^2(\Rb{+},L^2(Z,\nu,H))$ is possible thanks to the following result which is taken from
\cite[Theorem C.1]{Brz+Haus_2009}.
\begin{thm}
There exists a unique bounded linear operator
$$ I  : \mathcal{M}^2 (\mathbb{R}_+, L^2 (Z, \nu; H)) \rightarrow L^2 (\Omega,\mathcal{F}; H)$$
such that for $\xi \in  \mathcal{M}_{step}^2(\mathbb{R}_+, L^2(Z, \nu; H))$ we have $I(\xi ) =  \tilde{I}(\xi)$. In particular, there exists a constant
$C=C(H)$ such that for any $\xi\in \mathcal{M}^2(\Rb{+},L^2(Z,\nu,H))$,
\begin{equation}
\label{ineq-2.03} \mathbb{E} \vert \int_0^t \int_Z \xi(r,z)
\tilde{\eta}(dz,dr)\vert ^2_H \leq C\, \mathbb{E}
\int_0^t\int_Z\vert \xi(r,z)\vert^2_H\, \nu(dz)\,dr, \; t\geq 0.
\end{equation}
Moreover, for each $\xi\in \mathcal{M}^2(\Rb{+},L^2(Z,\nu,H))$ , the process $I (1_{[0,t]} \xi)$, $t \ge 0$, is an
$H$-valued \cadlag martingale. The process $1_{[0,t]}\xi$ is defined by $[1_{[0,t]} \xi ](r, z, \omega)
:= 1_{[0,t]} (r)\xi(r, z, \omega),$ $t\ge 0$, $r \in \mathbb{R}_+$ , $z\in Z$ and $\omega\in \Omega$.

As usual we will write
$$
 \int_0^t \int_Z \xi(r,z)
\tilde{\eta}(dz,dr) := I(\xi)(t),\quad t\ge 0.
$$

\end{thm}}

We will rewrite \eqref{1} in the following equivalent form
\begin{equation}  \label{2}
\begin{cases}
d\bu+\left[\kappa_1\mathcal{A}\bu+2\kappa_0\mathcal{A}_p \bu+B(\bu,\bu)\right]dt=\int_Z \sigma(t,\bu,z)\wie(dz,dt), \\
\bu(0)=\bu_0.%
\end{cases}%
\end{equation}
Here the operator $\mathcal{A}$ is defined through the relation
\begin{equation*}
\langle \mathcal{A}\bu,\bv\rangle=a(\bu,\bv)=\int_\mathcal{O} \frac{\partial \mathcal{E}_{ij}(\bu)}{%
\partial x_k}\frac{\partial \mathcal{E}_{ij}(\bv)}{\partial x_k} dx, \,\, \bu\in
D(A), \bv\in V.
\end{equation*}
Here and after the summation over repeated indices is enforced.

Note that
{
\begin{equation*}
D(\mathcal{A})=\lbrace \bu\in V: \exists f\in H\subset V^\ast \text{ for which }
a(\bu,\bv)=(f,\bv), \forall \bv\in V\rbrace,
\end{equation*}}
$\mathcal{A}=\mathbf{P}\Delta^2$, where $\mathbf{P}$ is the orthogonal projection
defined on $L^2(\mathcal{O})$ onto $H$.
\begin{Rem}\label{RM-1}
It is shown in \cite{BELLOUT5} that there exist two positive constants $k_1$ and $k_2$ depending only on $\mathcal{O}$ such that
\begin{equation}\label{eq-bloom}
 k_1 ||\bu||^2_{2}\le \langle \mathcal{A}\bu, \bu\rangle \le k_2 ||\bu||^2_{2},
\end{equation}
for any $\bu\in \ve$.
Thanks to this we will just write $\lVert \bu\rVert^2_2$ in place of $\langle \mathcal{A}\bu, \bu\rangle$, $\bu\in V$.
Also, it is not difficult to see that $\mathcal{A}$ is symmetric. This fact together with \eqref{eq-bloom} yields that
 $\mathcal{A}$ is self-adjoint.
\end{Rem}
\noindent We also introduce a bilinear form $B(\bu,\bv):\mathcal{H}^1\times \mathcal{H}^1
\rightarrow \mathcal{H}^{-1}$ as follows:
\begin{equation*}
\langle B(\bu,\bv),\bw\rangle=b(\bu,\bv,\bw)=\int_\mathcal{O}\bu_i\frac{\partial \bv_j}{\partial x_i}\bw_j dx,\,\,
\bu,\bv,\bw\in \mathcal{H}^1,
\end{equation*}
where $b(.,.,.)$ is the well-known trilinear form used in the mathematical
analysis of Navier-Stokes equations (see for instance \cite{temam}). The
bilinear form $B(\cdot,\cdot)$ enjoys the following properties:

\begin{itemize}
\item for any $\bu,\bv,\bw\in \mathcal{H}^1$, we have
\begin{equation}  \label{B1}
\langle B(\bu,\bv),\bw\rangle=-\langle B(\bu,\bw),\bv\rangle  \text{ and } \langle B(\bu,\bv),\bv\rangle=0.
\end{equation}

\item There exists a constant $C_0$ such that
\begin{equation}  \label{B2}
\langle B(\bu,\bv),\bw\rangle \le C_0 |\bu| ||\bv||_1||\bv||_2
\end{equation}
for any $\bu\in \mathcal{H}^1,\bv\in V, \bw \in V$.
\end{itemize}

The inequality \eqref{B2} can be proved by using H\"older's and Sobolev
inequalities (see \cite{temam}).

\noindent The nonlinear term $\mathcal{A}_p:V\rightarrow V^\ast$ is defined as
follows:
\begin{equation*}
(\mathcal{A}_p \bu, \bv)=\int_\mathcal{O}\Gamma(\bu)\mathcal{E}_{ij}(\bu)\mathcal{E}_{ij}(\bv) dx, \bu,\bv\in
V,
\end{equation*}
where $\Gamma(\bu)=(\varkappa+|\mathcal{E}(\bu)|^2)^{\frac{p-2}{2}}$. Some of
the properties of $\mathcal{A}_p$ is given below.

\begin{lem}\label{PROP-AP}
\begin{enumerate}
\item[(i)] There exists a positive constant $C(\varkappa, p)$ such that
\begin{equation}  \label{3}
||\mathcal{A}_p \bu-\mathcal{A}_p \bv||_{V^\ast}\le C ||\bu-\bv||_1,
\,\, \bu,\bv\in V.
\end{equation}

%\item[(ii)] If $\bu\in L^2(0,T,V)\cap L^\infty(0,T,H)$, then $\mathcal{A}_p \bu\in
%L^2(0,T,V^\ast)$.

\item[(ii)] For any $\bu,\bv\in V$
\begin{equation}  \label{AP3}
\langle \mathcal{A}_p \bu-\mathcal{A}_p\bv,\bu-\bv\rangle \ge 0.
\end{equation}
\end{enumerate}
\end{lem}
To check the results in the above lemma we need to recall the following results whose proofs can be found in
\cite{RUZICKA}.
\begin{lem}[\textbf{Korn's inequalities}]
 Let $1<p<\infty$ and let $\mathcal{O}\subset \mathbb{R}^d$ be of class $C^1$. Then there exist two positive constants $k^i_p=k^i(\mathcal{O},p), i=1,2$ such that
\begin{equation*}
 k_p^1 ||\bu||_{1}\le \left(\int_\mathcal{O} |\mathcal{E}(\bu)|^2 dx\right)^\frac{1}{2}\le k_p^2 ||\bu||_{1},
\end{equation*}
for any $\bu\in \mathcal{H}^1$.
\end{lem}
\begin{proof}[Proof of Lemma \ref{PROP-AP}]
Let $\bw$ be any element of $V$. Let us set $\delta=\frac{|\mathbf{e}|}{\sqrt{\varkappa}}$  and $\tilde{\varkappa}=\max(\varkappa^\frac{p-2}{2}, \varkappa^\frac{p-5}{2})$.
 Let us first note that
\begin{equation*}
\left| \frac{\partial \mathbf{T} }{\partial e_{ij}}(\mathbf{e})\right|\le \varkappa^\frac{p-2}{2} (1+\delta^2)^{\frac{p-2}{2}}+ \varkappa^\frac{p-5}{2}
\delta (1+\delta^2)^{\frac{p-4}{2}},
\end{equation*}
which implies that
\begin{equation*}%\label{PROP-AP1}
\left| \frac{\partial \mathbf{T} }{\partial e_{ij}}(\mathbf{e})\right|\le 3 \tilde{\varkappa}(1+\delta^2)^{\frac{p-2}{2}}.
\end{equation*}
Since $p\in (1,2]$ we have that
\begin{equation}\label{PROP-AP1}
\left| \frac{\partial \mathbf{T} }{\partial e_{ij}}(\mathbf{e})\right|\le 3 \tilde{\varkappa}.
\end{equation}
Secondly, we have
\begin{align*}
 \langle\mathbf{T}(\mathcal{E}(\bu))-\mathbf{T}(\mathcal{E}(\bv)), \bw\rangle&=2\kappa_0 \int_\mathcal{Q} \left[\mathbf{T}(\mathcal{E}(\bu))-\mathbf{T}(\mathcal{E}(\bv))\right]\cdot
 \mathcal{E}(\bw) dx,\\
&=2\kappa_0 \int_\mathcal{O} \mathcal{E}(\bw) \cdot \int_0^ 1 \frac{\partial \mathbf{T}(\mathcal{E}(\bv)+s(\mathcal{E}(\bu)-\mathcal{E}(\bv))) }{\partial s} ds dx ,\\
&=2\kappa_0 \int_\mathcal{O}  \mathcal{E}(\bw) \cdot (\mathcal{E}(\bu)-\mathcal{E}(\bv)) \int_0^1 \frac{\partial \mathbf{T}(\mathcal{E}(\bv)+s(\mathcal{E}(\bu)-\mathcal{E}(\bv))) }{\partial e_{ij}} ds dx .
\end{align*}
By using \eqref{PROP-AP1}  in the last equation yields
\begin{equation*}
 |\langle\mathbf{T}(\mathcal{E}(\bu))-\mathbf{T}(\mathcal{E}(\bv)), \bw\rangle|\le 6 \kappa_0 \int_\mathcal{O} \lvert \mathcal{E}(\bw)\rvert \cdot \lvert \mathcal{E}(\bu)-\mathcal{E}(\bv)\rvert dx.
\end{equation*}
Invoking H\"older's and Korn's inequalities implies the existence of a positive constant $K$ such that
\begin{equation*}
 |\langle\mathbf{T}(\mathcal{E}(\bu))-\mathbf{T}(\mathcal{E}(\bv)), \bw\rangle|\le 6 \tilde{\varkappa}\kappa_0 K ||\bu-\bv||_1 ||\bw||_1,
\end{equation*}
for any $\bu, \bv, \bw\in V$. We easily conclude from this the proof of (i).

It is known from \cite{ROKYTA} that for any $p\in(1,\infty)$
and for all $\mathbf{D}, \mathbf{E}\in
\mathbb{R}^{d\times d}_{sym}$:
\begin{align*}
\left(\mathbf{T}(\mathbf{D})-\mathbf{T}(\mathbf{E})\right)\cdot
\left(\mathbf{D}-\mathbf{E}\right)\ge 0.% \label{l4},
\end{align*}
Here
$$\mathbb{R}^{d\times d}_{sym}=\{\mathbf{D}\in \mathbb{R}^{d\times d}: D_{ij}=D_{ji}, i,j=1,2,\dots, d\}.$$
 Therefore, we see that for any $\bu,\mathbf{v}\in V$
\begin{align}
 \langle\mathcal{A}_p \bu-\mathcal{A}_p \mathbf{v}, \bu-\mathbf{v}\rangle=&\int_\mathcal{Q}\left[\mathbf{T}(\mathcal{E}(\bu))-\mathbf{T}(\mathcal{E}(\mathbf{v}))\right]\cdot \left[\mathcal{E}(\bu)-\mathcal{E}
(\mathbf{v})\right]dx\ge 0,\nonumber
\end{align}
which proves (ii).

\end{proof}
To close this section we introduce the main set of hypotheses used in this article. Throughout this work we suppose that we are given a function $\sigma$ satisfying the following set of constraints:
\begin{condition}\label{cond-s-1}
There  exist
%an intermediate space $\CV$,
nonnegative constants  %$\gamma$,
  $\ell_0,\ell_1, \ell_2$, $\ell_3$ such that, for any $t\in[0,T]$
and all $\bu,\bv\in H$, we have
\begin{enumerate}
%\item $|x|_{\CV} \le K\,|x | ^{1- \gamma}\Vert x \Vert ^{ \gamma}$ for all $x\in V$. % and $\gamma<\frac 12$.
%
\item  $|\sigma(t,\bu)|_{L^2(Z,\nu;H)}^2\le \ell_0+\ell_1|\bu|^2 $;
%\item { $|\sigma(t,u)|_{L^2(Z,\nu;H)}^2\le K_0+ K_1|u|^2_H +K_2|u|_{\CV}^2$};
%
\item  $|\sigma(t,\bu)-\sigma(t,\bv)|_{L^2(Z,\nu;H)}^2\le \ell_2|\bu-\bv|^2$.
%\item { $|\sigma(t,u)-\sigma(t,v)|_{L^2(Z,\nu;H)}^2\le L_1|u-v|^2_H +L_2|u-v|_{\CV}^2$}.
\item $|\sigma(t,\bu)|_{L^4(Z,\nu;H)}^4\le \ell_3(1+ |\bu|^4)$.
\end{enumerate}
\end{condition}

% The points (i) and (ii) were proved in \cite{GUO1}, and (iii) in \cite{DU}
% for example.  From now on, we will work with \eqref{2}.

\section{Existence of a strong solution}
This section is devoted to show that \eqref{1} admits at least one strong solution. The proof is based on Galerkin approximation and
idea borrowed from \cite{Breckner}. But before we proceed further we define {{explicitly}} what we mean by  strong solution of \eqref{1} or \eqref{2}.
\begin{Def}\label{mart-sol}
 Let $(Z,\CZ)$ be a separable metric space on which is defined a $\sigma$-finite measure $\nu$ and $\bu_0\in H$.
 A strong solution to the problem \eqref{2} is a stochastic process $\bu$ such that

\begin{enumerate}
  \item $\bu=\{\bu(t)\}_{t\ge 0}$ is a
$\mathbb{F}$-progressively measurable process such that
$$\mathbb{E}\sup_{s\in [0,T]}|\bu(s)|^4+\mathbb{E}\int_0^T ||\bu(t)||^2_2 dt <\infty,$$
\item the
following holds
\begin{equation}\label{weak-form}
 \begin{split}
 (\bu(t),\bw)=(\bu_0,\bw)-\kappa_1\int_0^t
\left( \langle\mathcal{A}\bu(s), \bw\rangle -\langle B(\bu(s),\bu(s)),\bw\rangle \right)ds\\-2\kappa_0\int_0^t \langle \mathcal{A}_p\bu(s),\bw\rangle ds+ \int_0^t\int_Z(\sigma(s,\bu(s),z),\bw) \wie(dz,ds),\\
 \end{split}
\end{equation}
for any $\bw\in V$, for almost all $t\in [0,T]$ and
$\PP$-almost surely.
\end{enumerate}
\end{Def}

\begin{thm}\label{main-ex}
Let the set of constraints in Condition \ref{cond-s-1} be satisfied. Then for any $T>0$ and %there exists a $K>0$ such that
initial value $\xi$ with $\EE|\xi|^4<\infty$,
there exists a solution $\bu=\{ \bu(t):0\le t<\infty\}$ to problem \eqref{2} which satisfies
%such that %for any $2\le p<2/\gamma $ and
\begin{equation}%\label{main-ex-1}
\EE \lk( \sup_{0\le t\le T} |\bu(t)| ^ {2r} + \int_0 ^ T \lk\Vert
\bu(s)\rk\Vert^2_{2} \lk| \bu(s) \rk| ^ {2r-2} \, ds \rk) \le C \lk(
\EE\lk|\xi\rk| ^ {2r} +1\rk), T\ge0,
\end{equation}
with $r=1,2$.
\end{thm}
Before we prove this result let us recall an important statement which is borrowed from \cite{Millet+Chueshov_2010}.

\begin{lem} \label{millet-lemma}
Let $X$, $Y$, $I$ and $\phi$ be non-negative processes and $Z$ be
a non-negative integrable random variable. Assume that $I$ is
non-decreasing and that there exist non-negative constants $C$,
$\alpha$,$\beta$, $\gamma$,  $\delta$ and $T$  satisfying first %with the following properties for some $T>0$
$$
\int_0 ^ T \phi(s)\, ds \le C,\quad \mbox{a.s.},\quad 2\beta e ^ C\le 1,\quad 2\delta e ^ C \le \alpha,
$$
and secondly for all $t\in[0,T]$ %$0\le t\le T$
there exists a constant $\tilde C>0$ such that \DEQS X(t) +\alpha
Y(t) \le Z +\int_0 ^ t\phi(r)X(r)\, dr +I(t),\quad  \mbox{a.s.},
\\
\EE I(t) \le \beta \EE X(t) +\gamma \int_0 ^ t \EE X(s)\, ds +\delta \EE Y(t) +\tilde C.
\EEQS
If $X\in L ^ \infty([0,T]\times \Omega)$, then we have
$$
\EE \lk[ X(t) +\alpha Y(t)\rk] \le 2\exp\lk( C+2t\gamma e ^ C\rk) \lk( \EE Z +\tilde C\rk),\quad t\in[0,T].
$$
\end{lem}

The proof of Theorem \ref{main-ex} will be split  into five steps.
\\[0.3cm]
%\paragraph{\bf 0 Step}
\textbf{ A priori uniform estimates:}\\
%\dela{
%Let $\{\phi_n:n\in\NN\}$ be an ONB in $H$ such that $\phi_n\in D(A)$, $H_n:= \span(\{\phi_i: i=1,\ldots n\}$,
%and $P_n:H\to H_n$ the orthogonal projection onto this basis.}
The operator $\mathcal{A}$ is self-adjoint and it follows from Rellich's theorem that
it is compact on $H$. Therefore, there exists a sequence of positive numbers
$\{\tilde{\lambda}_i: i = 1, 2, 3\dots \}$ and a family of smooth function $\{\phi_i : i = 1, 2, 3, \dots \}$ satisfying
\begin{equation}\label{eigen}
\mathcal{A}\phi_i  = \tilde{\lambda}_i \phi_i,
\end{equation}
for any $i\in \mathbb{N}$. We can assume that the family $\{\phi_i: i = 1, 2, 3, \dots \}$ so that it will
form an orthonormal basis of $H$  which is orthogonal and dense in $V$.

Let $\Pi_m$ denote the projection of $V^\ast$ onto $H_m:=span\{\phi_1,\cdots,\phi_m\}$. That is
$$\Pi_m x=\sum_{i=1}^m\langle x,\phi_i\rangle \phi_i, x\in \ve^\ast.$$
Also, ${\Pi_m}_{|_H}$ is the orthogonal projection of $H$ onto $H_m$.

For every $m\in\NN$, we consider the finite dimensional system of SDEs on $H_m$ given by
%\dela{
%\DEQSZ\label{discrete-sys}
%d (u_n(t) , \phi_k ) =  \la F( u_n(t)), \phi_k\ra \, dt +  \int_Z \la \sigma (u_n(t);z),\phi_k\ra \,\tilde \eta(dz,dt)
%,\quad k=1,\ldots n. %, \quad v\in H_n.
%\EEQSZ}

%\coma{
\DEQSZ\label{discrete-sys}
d \bum(t)  &=&  \Pi_m F( \bum(t))\, dt +  \int_Z \Pi_m\sigma (t,\bum(t),z) \,\tilde \eta(dz,dt),\ t\geq0,
\\ \bum(0)&=&\Pi_m\xi,\nonumber
\EEQSZ
where $F(\bum(s))=-A\bum(s)-\mathcal{A}_p \bum(s) + B(\bum(s), \bum(s))$. To shorten notation we set $B_m(\cdot, \cdot)=\prm B(\cdot, \cdot)$
 and $\sigma_m(\cdot, \cdot, \cdot)=\prm \sigma(\cdot, \cdot, \cdot)$.

We note that since $\Pi_m$ is a contraction of $V^\ast$, we
infer from \eqref{B1}, \eqref{B2} and point (2) of Condition \ref{cond-s-1} that $F$ is locally Lipschitz and
$\sigma_m:=\Pi_m\sigma$ is globally Lipschitz. As we know from
e.g. Albeverio, Brze{\'z}niak and Wu \cite{Alb+Brz+Wu_2010}, on
the basis of Condition \ref{cond-s-1},  equation
\eqref{discrete-sys} has a unique $H_m$-valued c\`{a}dl\`{a}g local  strong
solution $\bum$.  The following
proposition implies that it is in fact a global solution.
%}

\begin{prop}\label{proposition-estimate}
Let the assumptions be as in Theorem \ref{main-ex}. Then there exists a constant $C>0$ such that for $r=1,2$ we have
$$
\sup_m \EE \lk( \sup_{t\in[0,T]} |\bum(t)| ^ {2r} + \int_0 ^ T
\lk\Vert  \bum(s)\rk\Vert_2^2 \lk| \bum(s) \rk| ^ {2r-2} \, ds \rk)
\le C \lk( \EE\lk|\xi\rk| ^ {2r} +1\rk).
$$
\end{prop}
\begin{proof}[Proof of Proposition \ref{proposition-estimate}]
As mentioned above, it follows from  in \cite[Theorem 2.8]{Alb+Brz+Wu_2010} that Equation \eqref{discrete-sys} has a unique
c\`{a}dl\`{a}g local  strong solution $\bum$ in $H_m$. That
means that for any $m\in\NN$ there exists a unique solution on a short interval
interval $[0,T_m]$ satisfying
 \begin{equation*}
      \bum(t)  =\Pi_m\xi+\int_0^t \Pi_m F( \bum(s))\, ds +  \int_0^t\int_Z \sigma_m (s,\bum(s),z) \,\tilde \eta(dz,ds),\
 \quad t\in[0,T_m].
\end{equation*}
% Let us first of all assume that $\bum$ is bounded $\mathbb{P}$-a.s.
%
\medskip

We begin by checking the estimate in the proposition
with the case $r=1$.  { We argue as in \cite[Proof of
Theorem 3.1]{Alb+Brz+Wu_2010}.
Let $(\tau_M)_M$ be an incresing sequence of stopping times defined by   $$\tau_M=\inf\{t\in [0,T]: |\bum(t)|^2+ \int_0^t \lVert \bum(s)\rVert^2 ds\ge M^2 \}\wedge T,$$
for any integer $M$. Throughout we fix $r\in [0, T]$ and $0\le t\le r\wedge \tau_M$.} Since we can identify the space $H_m$ with $\mathbb{R}^m$ then
we can apply the finite dimensional It\^o's formula (see, for example, \cite[Chapter II, Theorem 5.1]{Ik-Wat-81})
 to the function $\vert \cdot\vert^{2r}$ and the process $\bum$. This procedure along with \eqref{B1} yields
\begin{equation}\label{prop-e-711}
\begin{split}
|\bum (t)|^2 = |\Pi_m\xi|^2 -
 2\kappa_1 \int_0^{t}\Vert \bum(s)\Vert^2_2\, ds -2 \kappa_0 \int_0^t \langle \mathcal{A}_p\bum(s), \bum(s) \rangle ds \\ +  2\int_0^{t}\int_Z ( \bum(s-),\sigma_m(s,\bum(s),z)) \tilde \eta(dz,ds)
 \\+ \int_0^{t}\int_Z \Psi(s,z)\eta(dz,ds),
 \end{split}
 \end{equation}
 where $$\Psi(s,z)=|\bum(s-)+\sigma_m(s,\bum(s),z)|^2-  |\bum(s-)|^2-\left( \bum(s-),
 \sigma(s,\bum(s),z)\right).$$
From the fact that $|x|^2 - |y|^2 + |x-y|^2 =2\la x-y,x\ra$, $x,y\in H$ and \eqref{AP3}, we derive from \eqref{prop-e-711} that
\begin{equation}\label{eq17}
\begin{split}
|\bum(t)|^2 + 2\int_0^{t}\kappa_1\Vert \bum(s)\Vert_2^2 ds \le
 |\Pi_m\xi|^2 + 2\int_0^{t}\int_Z ( \bum(s-),\sigma(s,\bum(s),z))
\tilde \eta(dz,ds)\\
+ \int_0^{t}\int_Z |\sigma(s,\bum(s),z)|^2\eta(dz,ds)
 \end{split}
 \end{equation}
 %\EEQS
 for any $t\in [0,r\wedge \tau_M]$, $r\in [0,T]$.
%  Thanks to \eqref{AP3} the above equation can be rewritten in the following form
%
%  \begin{equation}\label{eq17}
%   \begin{split}
%   |\bum(t)|^2 + 2 \kappa_1\int_0^{t}\Vert \bum(s)\Vert_2^2 ds  =
%   |\Pi_m\xi|^2 + { \int_0^{t}\int_Z |\sigma(s,\bum(s),z)|^2\eta(dz,ds)}\\
% + 2\int_0^{t}\int_Z ( \bum(s-),\sigma(s,\bum(s),z))\,
% \tilde \eta(dz,ds),
%  \end{split}
%  \end{equation}
% % \EEQS
%  for any $t \in [0,T]$.
For any $r\in [0, T]$ and $ t\in [0, r\wedge \tau_m]$ we define the following stochastic processes
\begin{align*}
X(t) :=&\sup_{0\le s\le t } |\bum(s)|^2;
\\
Y(t):=&\int_0^{t}\Vert \bum(s)\Vert_2^2\, ds;
\\
I(t) :=&\sup_{0\le s\le t} \Big( 2\Big| \int_0^{s}\int_Z (
\bum(\tau-),\sigma(\tau,\bum(\tau),z)) \tilde \eta(dz,d\tau)\Big|
+ \int_0^{s}\int_Z |\sigma(\tau,\bum(\tau),z)|^2\eta(dz,d\tau)\Big),\\
:=&\sup_{s\in [0,t]}|I_1(s)| + I_2(t),
 \end{align*}
 where
\begin{equation*}
 I_1(t)=\int_0^t\int_Z(\sigma(s,\buns,z),\buns)\wie(dz,ds),
\end{equation*}
and
\begin{equation*}
 I_2(t)=\sup_{0\le s\le t}\int_0^s \int_Z\vert\sigma(\tau,\bum(\tau-),z)\vert^2\eta(dz,d\tau).
\end{equation*}
Since $I_1(t)$ is a local martingale we can apply Burkholder-Davis-Gundy's inequality
 and get
\begin{equation*}
 \bae \sup_{s\in [0,r\wedge \tau_M]}|I_1(s)|\le C \bae\left(\int_0^{r\wedge \tau_m} \int_Z(\buns, \sigma(s,\bums,z))^2 \nu(dz)ds\right)^\frac{1}{2}.
\end{equation*}
    Thanks to H\"older's and Young's inequalities we have
\begin{align*}
  \bae \sup_{s\in [0,t]}|I_1(s)| &\le C \left[\eps \bae \sup_{s\in [0,t]}|\bums|^2 \right]^\frac{1}{2}\left[\frac{1}{\eps}\bae \int_0^t\int_Z|\sigma(s,\bums,z)|^2\nu(dz)ds \right]^\frac{1}{2}\nonumber\\
&\le C \eps \bae\sup_{s\in [0,t]}|\bums|^2+\frac{C}{\eps}\bae\int_0^t \int_Z |\sigma(s,\bums,z)|^2 \nu(dz)ds.
\end{align*}%
Invoking item (2) of Condition \ref{cond-s-1} we see that
\begin{equation}\label{eq19}
\bae \sup_{s\in [0,t]}|I_1(s)| \le C\eps \bae X(t) + \frac{C}{\eps} \ell_0 t +\frac{C}{\eps} \ell_1 \int_0^t \bae X(s) ds.
\end{equation}
Next, we will deal with the second term of $I(t)$. Taking into account that the process $$\int_0^t\int_Z|\sigma_m(r,\bum(r))|^2\eta(dz,dr)$$ has only positive jumps, we obtain
\begin{equation*}
\bae I_2(t)\le \bae \int_0^t \int_Z \vert\sigma(s,\bum,z)\vert^2 \nu(dz)ds.
\end{equation*}
Thanks to the item (1) of Condition \ref{cond-s-1} we see that
\begin{equation}\label{eq21}
\bae I_2(t)\le  \ell_0 t+ \ell_1 \int_0^t \bae X(s) ds.
\end{equation}
Thanks to \eqref{eq17} along with \eqref{eq19} and \eqref{eq21} we apply Lemma \ref{millet-lemma} and derive that there exist a positive constant $C$ such that
\begin{equation*}
\bae \left[\sup_{0\le s \le t}|\bum(s)|^2+2\kappa_1\int_0^t \Vert \bum(s)\Vert_2 ds \right]\le C (\bae \vert \xi\vert^2+1),
\end{equation*}
for any $m\in \mathbb{N}$ and $t\in [0,r\wedge \tau_m]$, $r\in [0,T]$.
We have just shown that
\begin{equation}\label{st-time}
\bae \left[\sup_{0\le s \le t\wedge \tau_M}|\bum(s)|^2+2\kappa_1\int_0^t \Vert \bum(s)\Vert_2 ds \right]\le C (\bae \vert \xi\vert^2+1) \forall t\in [0,T],
\end{equation}
with which we can infer that $$\mathbb{P}(\tau_M< t)\le C M^{-2}, \forall t\in
[0,T], \forall M>0.$$
Hence, $\lim_{M\rightarrow \infty }\mathbb{P}(\tau_M <t) = 0$, for
all $t\in [0,T]$. That is, $\tau_M \rightarrow \infty$ in
probability. Therefore, there exists a subsequence $\tau_{M_k}$
such that $\tau_{M_k}\rightarrow \infty$, a.s. Since the sequence
$(\tau_M)_M$ is increasing, we infer that $\tau_{M_k}\nearrow
\infty$ a.s.. Now we use Fatou's lemma and pass to the limit in \eqref{st-time} and derive that
\begin{equation*}
\bae \left[\sup_{0\le s \le t}|\bum(s)|^2+2\kappa_1\int_0^t \Vert \bum(s)\Vert_2 ds \right]\le C (\bae \vert \xi\vert^2+1).
\end{equation*}
The proposition is then proved for $r=1$.
Thus, it remains to show that it is true for the case $r=2$. We again apply It\^o's formula to obtain
 \begin{equation}\label{eqxy}
 \begin{split}
 |\bum(t)|^{2r}=
 2r \int_0^t \int_Z |\buns|^{2(r-1)}(\buns, \sigma_m(s,\bums,z))\wie(dz,ds)\\
 +|\Pi_m \xi|^{2r}+ 2r \int_0^t |\bum(s)|^{2(r-1)}(F_m(\bum(s)), \bum(s)) ds\\
 +\int_0^t\int_Z \Phi(s,z) \eta(dz,ds),
 \end{split}
 \end{equation}
 where $$ \Phi(s,z)=|\buns+\sigma_m(s,\bums,z)|^{2r}-|\buns|^{2r}-2r |\buns|^{2(r-1)}(\buns, \sigma(s,\bums,z)). $$
 Thanks to \eqref{B1} and \eqref{AP3} the estimate \eqref{eqxy} becomes
 \begin{equation*}
 \begin{split}
 |\bum(t)|^{2r}+2r \kappa_1 \int_0^t \vert \bums|^{2(r-1)}\Vert \bums\Vert_2^2 ds -|\Pi_m \xi|^{2r}-\int_0^t\int_Z \Phi(s,z) \eta(dz,ds)\\
\le 2r \int_0^t \int_Z |\buns|^{2(r-1)}(\buns, \sigma_m(s,\bums,z))\wie(dz,ds).
 \end{split}
 \end{equation*}
 Taking the supremum over $[0,t]$ on both sides of the above estimate leads to
 \begin{equation}\label{eq25}
 \begin{split}
 \sup_{s\in [0,t]}|\bum(s)|+2r \kappa_1\int_0^t \vert \bums|^{2(r-1)}\Vert \bums\Vert_2^2 ds
 \le
 |\Pi_m \xi|^{2r}+ J(t),
 \end{split}
 \end{equation}
 where $J(t)=J_1(t)+J_2(t)$ with
 \begin{align*}
 J_1(t)&=2r \sup_{s\in [0,t]}\left\vert\int_0^s\int_Z |\buns|^{2(r-1)}(\bu(\tau-), \sigma_m(s,\bum(\tau),z))\wie(dz,d\tau)\right\vert,\\
 J_2(t)&=\sup_{s\in [0,t]}\left\vert\int_0^s\int_Z \Phi(\tau,z) \eta(dz,d\tau)\right\vert.
 \end{align*}
 First, we apply the  Burkh{o}lder-Davis-Gundy  inequality
 \begin{eqnarray*}
    \mathbb{E}J_1(t)
\leq 2r C\mathbb{E}\left(\int_0^{t}\int_Z|\bum(s)|^{4(q-1)}|\bum(s)|^2 |\sigma_m(s,\bum(s),z))|^2 \,\nu(dz)ds\right)^{\frac{1}{2}}.
\end{eqnarray*}
 Then using item (1) of Condition  \ref{cond-s-1} and H\"older's
inequality implies.
\begin{equation*}
\begin{split}
\mathbb{E}J_1(t)
\leq \left(\frac{1}{\varepsilon}\EE\int_0^{t}|\bum(s)|^{2r-2}\Big{(}\ell_0+\ell_1|\bum(s)|^2\Big{)}ds\right)^{\frac{1}{2}}
\\ \times 2rC\left(\varepsilon\mathbb{E}\sup_{s\in[0,t]}|\bum(s)|^{2r}\right)^{\frac{1}{2}}.
\end{split}
\end{equation*}
Invoking Young's inequality yields
\begin{equation*}%\label{eq28}
\begin{split}
\mathbb{E}J_1(t)
\leq \frac{2rC}{2}\varepsilon\mathbb{E}\sup_{s\in[0, t]}|\bum(s)|^{2r}+\frac{r C \ell_0}{2\varepsilon}\mathbb{E}\int_0^{t}|\bum(s)|^{2r-2}ds+\frac{2rC \ell_1}{2\varepsilon}\mathbb{E}
\int_0^{t}|\bum(s)|^{2r}ds.
\end{split}
\end{equation*}
Using the fact that
for $r\ge 2$, $|x|^{2r-2}\leq 1+|x|^{2r}$ we deduce from the last inequality that
\begin{equation}\label{eqkn0}
\begin{split}
    \mathbb{E}J_1(t)
\leq\frac{2rC \ell_0T}{2\varepsilon}+\frac{2rC}{2}\varepsilon\mathbb{E}\sup_{s\in[0, t]}|\bum(s)|^{2r}+\frac{2rC(\ell_0+\ell_1)}{2\varepsilon}\mathbb{E}\int_0^{t}|\bum(s)|^{2r}ds.
\end{split}
\end{equation}
Now we deal with $J_2(t)$. First, note that

\begin{equation}\label{eqkn}
     \mathbb{E}J_2(t)
\leq\left(\frac{r^2+r}{2}\right)\mathbb{E}\int_0^{t}\int_Z\Big{(}|\bum(s-)|^{2(r-1)}|\sigma_m(s,\bum(s),z)|^2
+|\sigma_m(s,\bum(s),z)|^{2r}\Big)\nu(dz)ds,
\end{equation}
where we have used the fact that
 \begin{equation}\label{taylor-1}
     \left||x+h|^{2r}-|x|^{2r}-2r|x|^{2(r-1)}(x,h)\right|\leq \frac{r^2+r}{2}(|x|^{2(r-1)}|h|^2+|h|^{2r}),
\end{equation}
\text{for all }$x,h\in H$. Let set $C_r=\frac{r^2+r}{2}$. Now thanks to items (1) and (3) of Condition \ref{cond-s-1} we derive from \eqref{eqkn} that there exist positive constants
$\ell_r$ and $C_r$ such that
\begin{eqnarray}
\bae J_2(t) &\leq&C_r \ell_r\mathbb{E}\int_0^{t}(1+|\bum(s)|^{2r}+|\bum(s)|^{2(r-1)}\Vert \bum(s)\Vert_2 ^{2r})ds\nonumber\\
&\leq&C_r \ell_r T+C_r \ell_r\int_0^{t}|\bum(s)|^{2r}ds+C_r\ell_r\int_0^t|\bum(s)|^{2(r-1)}\Vert
\bum(s)\Vert_2 ^2ds.\label{eqkn2}
\end{eqnarray}
Therefore we see from \eqref{eqkn0} and \eqref{eqkn2} that there exist positive constants $C_r^\prime$, $M_r$, $\ell_r^\prime$, and $L^\prime_r$ such that
\begin{equation}\label{eq33}
\bae J(t)\le C_r^\prime T + M_r \eps \bae \sup_{s\in[0,t]} |\bums|^{2r}+ \frac{\ell_r^\prime}{2 \eps}\int_0^t \bae |\bums|^{2r} ds + L_r^\prime \bae \int_0^t  |\bums|^{2(r-1)}\Vert \bums \Vert ^2_1 ds,
\end{equation}
for any $m\in \mathbb{N}$ and $t\in [0,T]$. Let set  $$X(t)=\sup_{s \in [0,t]} |\bums|^{2r},$$ and
$$Y(t)=\int_0^t |\bums|^{2(r-1)} \Vert \bums \Vert_2^2 ds.$$
Thanks to \eqref{eq25}, \eqref{eq33} and an appropriate choice of $\eps>0$ we find that $X(\cdot)$ and $Y(\cdot)$ verifies the conditions in Lemma \ref{millet-lemma}. Therefore we infer the existence of a positive constant $C$ such that
\begin{equation*}%\label{discrete-sys}
\bae \sup_{s \in [0,t]}|\bums|^{2r}+ C(q, \kappa_1) \int_0^t |\bums|^{2(r-1)} \Vert \bums\Vert_2^2ds \le C (\bae |\xi|^{2r} +1 ), r=2,
\end{equation*}
for any $t\in [0,T]$ and $m \in \mathbb{N}$.
This completes the proof of the proposition.
\end{proof}

\textbf{Passage to the limit:}

To prove the existence of the solution of \eqref{2} we need to pass to the limit in the terms of \eqref{discrete-sys} and in the estimate of Proposition \ref{proposition-estimate}.
Before we do so we recall that there exists a constant $C>0$ such that
\begin{align}
\sup_{m \in \NN}\left(\bae \sup_{s\in [0,T]} |\bums|^{2r}+ \bae \int_0^T \lvert \bums\rvert^{2(r-1)}\Vert \bums\Vert^2_2ds\right)< C.\label{LIM1}\\
%\sup_{m \in \NN}\left(\bae \int_0^T \vert \bums\vert \Vert \bums\Vert^2_2ds\right)< C.\label{LIM2}
\end{align}
We have the following weak compactness result.
\begin{prop}\label{WEAK-COMP}
We can extract from $\bum$ a subsequence which is not relabeled and there exists a stochastic process $\bu$ such that
\begin{align}
\bum &\rightharpoonup \bu \text{(weak star) in } L^4(\Omega; L^\infty([0,T]; H)),\label{LIM3}\\
\bum &\rightharpoonup \bu \text{ in } L^2(\Omega\times [0,T]; V).\label{LIM4}
\end{align}
Moreover, there exists three elements $\bb, \Sigma, \mathbf{A}$ such that
\begin{align}
B_m(\bum, \bum) &\rightharpoonup \bb \text{ in } L^2(\Omega\times [0,T]; V^\ast),\label{LIM5}\\
\mathcal{A}_p\bum &\rightharpoonup \mathbf{A} \text{ in } L^2(\Omega\times [0,T]; V^\ast),\label{LIM6}
\\
\sigma_m(t,\bum,\cdot) &\rightharpoonup \Sigma \text{ in } L^2(\Omega\times [0,T]; L^2(Z, \nu; H)).\label{LIM7}
\end{align}
\end{prop}
\begin{proof}
Since $L^2(\Omega\times [0,T]; V)$ is a Hilbert space and $L^4(\Omega; L^\infty(0,T; H))$ is a Banach space and the dual of $L^\frac{4}{3}(\Omega; L^1(0,T; H)$,
 we easily infer from Banach-Alaoglu's theorem  and the uniqueness of weak limit that there exist a subsequence of $\bum$
 (which is denoted with the same fashion) and a stochastic process $\bu$ belonging to $L^4(\Omega; L^\infty(0,T; H))\cap L^2(\Omega\times[0,T]; V)$
 such that \eqref{LIM3} and \eqref{LIM4} hold true.

It remains to show \eqref{LIM5}-\eqref{LIM7}. To prove \eqref{LIM5} we first recall that there exists a positive constant $C_0$ such that
\begin{equation*}
\vert\langle B_m(\Phi, \bv), \bw\rangle\vert\le C_0 \vert \Phi\vert \Vert \bv \Vert_1 \Vert \bw\Vert_2,
\end{equation*}
for any $\Phi\in H$, $\bv\in V$ and  $\bw\in V$.
This inequality implies that
\begin{equation*}%\label{LIM8}
\bae \int_0^T \Vert B_m(\bum(s), \bum(s))\Vert^2_{V^\ast}\le C \bae \int_0^T |\bums|^2 \Vert \bums \Vert_1^2 ds,
\end{equation*}
from which and \eqref{LIM1} we get that $B_m(\bum, \bum)$ is a bounded sequence in the Hilbert space $L^2(\Omega\times[0,T]; V^\ast)$.
Thus, there exists an element of $L^2(\Omega\times[0,T]; V^\ast)$ that
we denote by $\bb$ such that $B_m(\bum, \bum)$ converges weakly to $\bb$ in $L^2(\Omega\times[0,T]; V^\ast)$.

By invoking \eqref{3} and \eqref{LIM1} we see that the following uniform estimate holds
\begin{equation*}
\sup_{m\in \NN}\bae \int_0^T \Vert \mathcal{A}_p \bums\Vert^2_{V^\ast} ds \le C.
\end{equation*}
Therefore, the proof of \eqref{LIM6} follows the same lines as for the proof of \eqref{LIM5}.

From point (1) of Condition \ref{cond-s-1} and estimate \eqref{LIM1} we easily obtain the uniform estimate
\begin{equation*}
 \sup_{m\in \NN}\bae \int_0^T\Vert \sigma_m(s,\bums,z)\Vert^2_{L^2(Z,\nu;H)}\le K_0T +K_1\bae \int_0^T |\bums|^2 ds+K^\prime_1\bae \int_0^T \Vert \bums \Vert_2 ds\le C,
\end{equation*}
which implies that  $\sigma_m(s,\bums,z)$ is a bounded sequence in $L^2(\Omega\times [0,T]; L^2(Z,\nu;H))$. Therefore, by Banach-Alaoglu we deduce
the existence of $\Sigma$ belonging to $L^2(\Omega\times [0,T]; L^2(Z,\nu;H))$ such that \eqref{LIM7} holds. This completes the proof of the proposition.
\end{proof}
With the convergences in Proposition \ref{WEAK-COMP} we can pass to the limit in each term of \eqref{discrete-sys} and obtain that
\begin{equation}\label{LIM9}
 \bu(t)+\kappa_1 \int_0^t  \mathcal{A} \bu(s) ds+\kappa_0 \int_0^t \mathbf{A}(s) ds = \bu_0 + \int_0^t \int_Z \Sigma(s,z) \wie(dz,ds),
\end{equation}
$\PP$-a.s. and for any $t\in [0,T]$ as an equality in $V^\ast$. Also, passing to the limit in \eqref{LIM1} gives the estimate in Theorem \ref{main-ex}.
Thanks to \eqref{LIM5} and \eqref{LIM7} we can deduce from  \cite{Gyongy-Krylov} that the stochastic process $\bu$ has a \cadlag modification taking values in $H$. From now on we will
identify $\bu$ with its \cadlag modification. Henceforth,  we need to show the following identities to complete the proof of Theorem  \ref{main-ex}.
\begin{prop}\label{WEAK-COMP2}
 We have the following identities
\begin{align}
 \bb&=B(\bu,\bu) \text{ in } L^2(\Omega\times [0,T]; V^\ast),\label{LIM10}\\
\mathbf{A}&=\mathcal{A}_p \bu \text{ in } L^2(\Omega\times [0,T]; V^\ast),\label{LIM11}\\
\Sigma&= \sigma(t,\bu,\cdot) \text{ in } L^2(\Omega\times [0,T]; L^2(Z,\nu;H))\label{LIM12}.
\end{align}
\end{prop}
For any integer $M\ge 1$ we consider the sequence of stopping times $\{\tau_M: M\ge 1\}$ defined by
\begin{equation*}
 \tau_M= \inf\{ t \in [0,T]: |\bu(t)|^2+\int_0^T \Vert \bu(s)\Vert_2^2\ge M^2 \}\wedge T.
\end{equation*}
The proof of Proposition \ref{WEAK-COMP2} requires the following convergences.
\begin{lem}\label{COMP}
 For any $M\ge 1$ we have that, as $m\rightarrow \infty$,
\begin{equation}\label{LIM13}
 1_{[0,\tau_M]}\left( \bum-\bu\right)\rightarrow 0 \text{ in } L^2(\Omega\times [0,T]; V),
\end{equation}
and
\begin{equation}\label{LIM14}
 \bae |\bum(\tau_M) -\bu(\tau_M)|\rightarrow 0.
\end{equation}
\end{lem}
\begin{proof}[Proof of Lemma \ref{COMP}]
 Let $\bbu$ be the orthogonal projection of $\bu$ onto $\Span\{\phi_1, \dots, \phi_m\}$, that is $$\bbu=\sum_{i=1}^m (\bu,\phi_i)\phi_i.$$  It is clear that as $m\rightarrow \infty$
\begin{equation}\label{LIM15}
 \bbu\rightarrow \bu \text{ in } L^2(\Omega\times [0,T]; H).
\end{equation}
We also can check easily that
\begin{equation}\label{54-1}
 \bae |\bbu(\tau_M)-\bu(\tau_M)|^2\rightarrow 0,
\end{equation}
as $m\rightarrow \infty$.

First we should note that
\begin{align*}
 \langle \mathcal{A} \bbu(t), \bbu(t)\rangle=\langle \sum_{j} \mathcal{A}\phi_j (\bbu(t), \phi_j), \sum_{i} (\bbu(t), \phi_i)\phi_i\rangle
\end{align*}
Thanks to \eqref{eigen} we have
\begin{align*}
  \langle \mathcal{A} \bbu(t), \bbu(t)\rangle=&\sum_{i,j}(\bbu(t), \phi_j)(\bbu(t), \phi_i) (\lambda_j \phi_j, \phi_i),\\
  =& \sum_{j}(\bbu(t), \phi_j)^2(\lambda_j \phi_j, \phi_j).
\end{align*}
Thanks to \eqref{eigen} again we have
\begin{align*}
 \langle \mathcal{A} \bbu(t), \bbu(t)\rangle=\sum_{j}(\bbu(t), \phi_j)^2\langle \mathcal{A} \phi_j, \phi_j\rangle.
\end{align*}
From this, we can easily derive that
\begin{equation}\label{LIM16}
||\bbu(t)||^2_2\le ||\mathcal{A}||\,\, |\bu(t)|^2,
\end{equation}
for almost all $(\omega, t)\in \Omega \times [0,T].$
Also,
\begin{align}
 \Vert \bbu(s)-\bu(s)\Vert^2_2\le &\langle \mathcal{A}\bbu(s)-\mathcal{A}\bu(s), \bbu(s)-\bu(s)\rangle\nonumber\\
 \le & \langle \sum_{i=m+1}^\infty (\bu(s),\mathcal{A}\phi_i) \phi_i, \sum_{j=m+1}^\infty (\bu(s),\phi_j) \phi_j \rangle,\nonumber\\
\le & \langle\sum_{i=m+1}^\infty (\bu(s),\phi_i) \tilde{\lambda_i}\phi_i ,\sum_{j=m+1}^\infty (\bu(s),\phi_j)\phi_j\rangle,\nonumber\\
\le & \sum_{j=m+1}^\infty (\bu(s),\phi_j)^2 \langle \mathcal{A} \phi_j, \phi_j\rangle\nonumber\\
\le & ||\mathcal{A}||\sum_{j=m+1}^\infty (\bu(s),\phi_j)^2, \label{44-b}
\end{align}
for any $m$. Since $\bu\in H$ for almost all $(\omega, t)\in \Omega\times [0,T]$, we see that the right hand side of the last inequality converges to 0 as $m\rightarrow \infty$.
Therefore
 \begin{equation}\label{LIM16-a}
 \bbu(s)\rightarrow \bu \text{ in $V$ for almost all $(\omega, t)\in \Omega\times [0,T]$. }
 \end{equation}
  Furthermore, owing to \eqref{LIM16} and the dominated convergence theorem we can state that
 \begin{equation}\label{LIM17}
  \bbu\rightarrow \bu \text{ in } L^2(\Omega\times [0,T]; V).
 \end{equation}
Next, it is not difficult to see that $\bbu$ satisfies the following equations
\begin{equation*}%\label{LIM18}
 \bbu(t)+\kappa_1 \int_0^t \mathcal{A}\bbus ds+ \kappa_0 \int_0^t \prm \mathbf{A}(s) ds+\int_0^t \prm B(s)= \prm \xi+\int_0^t \int_Z \prm \Sigma(s,z) \wie(dz,ds).
\end{equation*}
Let $\xm$ be the stochastic processes defined by $\xm=\bum-\bbu$. From the equations of the last line and \eqref{discrete-sys} we obtain
\begin{equation*}
 \begin{split}
  \xm(t) +\kappa_0 \int_0^t [\mathcal{A}_p \bums-\prm \mathbf{A}(s)]ds+\int_0^t [B_m(\bums,\bums)-\prm \bb(s)]ds\\
= \int_0^t \int_Z [\sigma_m(s,\bums, z)-\prm \Sigma(s,z)]\wie(dz,ds)-\kappa_1\int_0^t \mathcal{A} \xms ds.
 \end{split}
\end{equation*}
 Applying It\^o's formula to the funxtion $\Phi(x)=|x|^2$ and $\xm(t)$ yields
\begin{equation*}
 \begin{split}
  |\xm(t)|^2+2\kappa_1\int_0^t \langle \mathcal{A} \xms, \xms\rangle ds+2\kappa_0\int_0^t \langle \mathcal{A}_p \bums-\prm \mathbf{A} (s), \xms\rangle ds\\
=2\int_0^t \langle \prm \bb(s)-B_m(\bums,\bums), \xms\rangle ds+\int_0^t\int_Z \Psi(s,z) \eta(dz,ds)\\
+2\int_0^t \int_Z (\sigma_m(s,\bums,z)-\prm \Sigma(s,z), \xms) \wie(dz,ds),
 \end{split}
\end{equation*}
where
\begin{align*}
\Psi(s,z)=&|\xm(s-)+\sigma_m(s,\bums,z)-\prm \Sigma(s,z)|^2-|\xm(s-)|^2\\
& \quad -2(\sigma_m(s,\bums,z)-\prm\Sigma(s,z), \xm(s-))\\
           =& |\sigma_m(s,\bums,z)-\prm \Sigma(s,z) |^2.
\end{align*}
Let $r(t)$ be the real valued stochastic process defined by $r(t)=K_1 t+\frac{C_0^2}{4 \kappa_1}\int_0^t ||\bu(s)||^2_2ds$.
 Applying It\^o's formula to $e^{-r(t)} |\xm(t)|^2$ leads to
\begin{equation}\label{LIM19}
 \begin{split}
    e^{-r(t)}|\xm(t)|^2+2\kappa_1\int_0^te^{-r(s)} \Vert\xms\Vert_2^2ds+2\kappa_0\int_0^te^{-r(s)} \langle \mathcal{A}_p \bums-\prm \mathbf{A} (s), \xms\rangle ds\\
=2\int_0^t e^{-r(s)}\langle \prm \bb(s)-B_m(\bums,\bums), \xms\rangle ds-\frac{C_0^2}{4 \kappa_1}\int_0^t e^{-r(s)}|\xms|^2 ||\xms||^2_2 ds\\
 - K_1\int_0^t e^{-r(s)} |\xms|^2 ds+\int_0^t\int_Z e^{-r(s)} |\sigma_m(s,\bums,z)-\prm \Sigma(s,z) |^2\eta(dz,ds)\\
+2\int_0^t e^{-r(s)}\int_Z (\sigma_m(s,\bums,z)-\prm \Sigma(s,z), \xms) \wie(dz,ds),
 \end{split}
\end{equation}
Let us study each term term of \eqref{LIM19}. For the nonlinear term involving $B_m$ and $\bb$ we have that
\begin{equation}
 B_m(\bbu, \bbu)-B_m(\bum, \bum)=B_m(\xm, \bbu)+B_m(\bum,\xm).\label{guy-14}
\end{equation}
Out of this and \eqref{B1} we obtain that
\begin{equation*}
\begin{split}
 \langle\prm \bb(s)-B_m(\bums,\bums), \xms \rangle=\langle\prm \bb(s)-B_m(\bbus, \bbus), \xms \rangle\\ +\langle B_m(\xms, \bbus),\xms \rangle,
\end{split}
\end{equation*}
which along with \eqref{B2} and Young's inequality  imply that
\begin{equation}\label{NON-B}
 \begin{split}
  \langle\prm \bb(s)-B_m(\bums,\bums), \xms \rangle\le \langle\prm \bb(s)-B_m(\bbus, \bbus), \xms \rangle\\
\frac{C^2_0}{4 \kappa_1} |\xms|^2||\bbus||^2_2 +\kappa_1 ||\xms||^2_2.
 \end{split}
\end{equation}
Next, we have
\begin{equation}\label{NON-A1}
\begin{split}
 \langle \mathcal{A}_p \bums-\mathbf{A}(s), \xms\rangle =\langle \mathcal{A}_p \bums-\mathcal{A}_p \bbus, \xms \rangle\\
+ \langle\mathcal{A}_p \bbus-\mathbf{A}(s), \xms \rangle.
\end{split}
\end{equation}
Invoking the point (ii) of Lemma \ref{PROP-AP} we see  that
\begin{equation}\label{NON-A2}
 \langle \mathcal{A}_p \bums-\mathcal{A}_p \bbus,\xms\rangle\ge 0.
\end{equation}
Setting $S=|\sigma_m(s,\bums,z)-\prm \Sigma(s,z)|^2$ we see that
\begin{equation*}
 \begin{split}
  S=|\prm [\sigma(s,\bums,z)-\sigma(s,\bu(s),z)]|^2-|\prm [\sigma(s,\bu(s),z)-\Sigma(s,z)]|^2\\
2\left( \prm [\sigma(s,\bums,z)-\Sigma(s,z)], \prm [\sigma(s,\bu(s),z)-\Sigma(s,z)]\right).
 \end{split}
\end{equation*}
Owing to point (1) of Condition \ref{cond-s-1} we have that
\begin{equation}\label{NON-S}
 \begin{split}
  S\le \ell_2|\xms|^2+\ell_2 |\bbus-\bu(s)|^2- |\prm [\sigma(s,\bu(s),z)-\Sigma(s,z)]|^2\\
2\left( \prm [\sigma(s,\bums,z)-\Sigma(s,z)], \prm [\sigma(s,\bu(s),z)-\Sigma(s,z)]\right).
 \end{split}
\end{equation}
Putting \eqref{NON-B}, \eqref{NON-A1}, \eqref{NON-A2} and \eqref{NON-S} into \eqref{LIM19}, replacing $t$ by $\tau_M$ and taking the mathematical expectation lead to
\begin{equation}\label{LIM20}
 \begin{split}
   \bae e^{-r(\tau_M)}|\xm(\tau_M)|^2+\bae \int_0^{\tau_M} \int_Z e^{-r(s)} |\prm [\sigma(s,\bu(s),z)-\Sigma(s,z) ]|^2 \eta(dz,ds)
\\ \le -\kappa_1\bae \int_0^{\tau_M}e^{-r(s)} ||\xms||^2_2 ds
+ 2\kappa_0 \bae \int_0^{\tau_M} e^{-r(s)} \langle \prm \mathbf{A}(s) -\mathcal{A}_p \bbus, \xms \rangle ds\\
+ 2\bae\int_0^{\tau_M} \int_Ze^{-r(s)}  \left( \prm [\sigma(s,\bums,z)-\Sigma(s,z)], \prm [\sigma(s,\bu(s),z)-\Sigma(s,z)]\right)\eta(dz,ds)\\
+
\bae \int_0^{\tau_M} e^{-r(s)}
\langle \prm \bb(s)-B_m(\bbus,\bbus), \xms\rangle ds\\
+K_1 \bae \int_0^{\tau_M} |\bbus-\bu(s)|^2 e^{-r(s)} ds.
 \end{split}
\end{equation}
Now we will show that the last four terms of the right hand side of \eqref{LIM20} will tend to 0 as $m\rightarrow 0$.
Thanks to \eqref{LIM15} we have
\begin{equation}\label{LIM21}
 \bae \int_0^T 1_{[0,\tau_M]}(s) e^{-r(s)} |\bbus-\bu(s)|^2 ds \rightarrow 0.
\end{equation}
Owing to \eqref{guy-14} and \eqref{B2} we see that
\begin{equation}\label{53-b}
\begin{split}
\left \Vert 1_{[0,\tau _{M}]}(t) e^{-r(t)} [B(\bbu(t),\bbu(t))-%
B(\bu(t),\bu(t))]\right\Vert _{V^{\ast }}\leq & 1_{[0,\tau
_{M}]}(t)C_{0}||\bbu(t)||_1\,\,|\bbu(t)-\bu(t)| \\
& +1_{[0,\tau _{M}]}(t)C_0 |\bu(t)| \,\, ||\bbu(t)-\bu(t) ||_2,
\end{split}%
\end{equation}%
which with \eqref{LIM16-a} implies that
\begin{equation*}
\left\Vert 1_{[0,\tau _{M}]}(t)e^{-r(t)}[B(\bbu(t),\bbu(t))-%
B(\bu(t),\bu(t))]\right\Vert _{V^{\ast }}\rightarrow 0\text{
a.e. }(\omega ,t)\in \Omega \times \lbrack 0,T],  %\label{CONV8}
\end{equation*}%
as $m\rightarrow \infty $. Furthermore, owing to \eqref{LIM16} and \eqref{44-b} we see from \eqref{53-b} that
\begin{equation}
\left \Vert 1_{[0,\tau _{M}]}(t) e^{-r(t)} [B(\bbu(t),\bbu(t))-%
B(\bu(t),\bu(t))]\right\Vert _{V^{\ast }}\leq 2 C_0 M ||A||^\frac{1}{2}\lvert\bu(t)\rvert. %\in L^2(\Omega\times[0,T]; \mathbb{R})
\label{CONV9}
\end{equation}%
Note that $\lvert\bu(t)\rvert$ is bounded in $\in L^{2}(\Omega \times \lbrack 0,T],\mathbb{R})$. Thus,
the Dominated Convergence Theorem implies that
\begin{equation}
\left\Vert 1_{[0,\tau _{M}]}(t)e^{-r(t)}[B(\bbu(t),\bbu(t))-%
B(\bu(t),\bu(t))]\right\Vert _{V^{\ast }}\rightarrow 0\text{ in } L^2(\Omega\times [0,T]; \mathbb{R} )\label{CONV8-b}
\end{equation}%
 By the convergences \eqref{LIM4} and \eqref{LIM17} we have
\begin{equation}\label{WEAK}
\bbu-\bum \rightharpoonup 0 \text{ in } L^{2}(\Omega;L^{2}(0,T;V)).
\end{equation}%
We derive from this, (\ref{CONV9}) and (\ref{CONV8-b}) that
\begin{equation*}
E\int_{0}^{\tau _{M}}e^{-r(s)}\langle B(\bbus,\bbus)-B%
(\bu(s),\bu(s)),\bbus-\bums\rangle ds\rightarrow 0
\end{equation*}%
as $m\rightarrow \infty $. Hence
\begin{equation*}
\begin{split}
& \lim_{m\rightarrow \infty }E\int_{0}^{\tau _{M}}e^{-r(s)} \langle B%
(\bbus,\bbus)-B^{\ast }(s),\bbus-\bums\rangle ds \\
=& \lim_{m\rightarrow \infty }E\int_{0}^{\tau _{M}}e^{-r(s)}\langle B%
(\bbus,\bbus)-B(\bu(s),\bu(s)),\bbus-\bums\rangle ds \\
& +\lim_{m\rightarrow \infty }E\int_{0}^{\tau _{M}}e^{-r(s)}\langle B%
(\bu,\bu)-B^{\ast }(s),\bbus-\bums\rangle ds \\
=& 0.
\end{split}
%\label{CONV10}
\end{equation*}%
Since $\prm\circ \prm =\prm$ and $\lVert \prm \rVert\le 1 $, it follows that $1_{[0,\tau_M]}e^{-r(s)}\prm[\sigma(s,\bu(s),z)-\Sigma(s,z)]$
is bounded in $L^2(\Omega\times [0,T]; L^2(Z,\nu;H))$. Therefore we see from \eqref{LIM7} that
\begin{equation*}%\label{CONV11}
 2\bae\int_0^{\tau_M} \int_Ze^{-r(s)}  \left( \prm [\sigma(s,\bums,z)-\Sigma(s,z)], \prm [\sigma(s,\bu(s),z)-\Sigma(s,z)]\right)\eta(dz,ds) \rightarrow 0
\end{equation*}
as $m\rightarrow \infty$.

Now it is not difficult to check that
\begin{equation*}
\begin{split}
 \bae \int_0^t e^{-r(s)} \langle\prm [\mathbf{A}(s)-\mathcal{A}_p\bbus], \xms \rangle ds=\bae \int_0^t e^{-r(s)} \langle\prm [\mathbf{A}(s)-\mathcal{A}_p \bu(s)], \xms \rangle ds\\+
\bae \int_0^t e^{-r(s)} \langle\prm [\mathcal{A}_p\bu(s)-\mathcal{A}\bbus], \xms \rangle ds.
\end{split}
\end{equation*}
Since $\langle \prm \bv, \xm\rangle=\langle \bv, \xm\rangle$ for any $\bv \in V^\ast$   and $e^{-r(s)} (\mathbf{A}(s)-\mathcal{A}_p \bu(s))$ is a bounded element of $L^2(\Omega\times[0,T]; V^\ast)$ , we derive from \eqref{WEAK}
that the first term of the right hand side of the  above equation tends to zero as $m\rightarrow \infty$. Owing to point (i) of Lemma \ref{PROP-AP}, the strong convergence
 \eqref{LIM17} and the weak convergence \eqref{WEAK} we see that the second term of the right hand side converges to zero as well. Thus, we have just proved that
\begin{equation*}%\label{CONV12}
 \bae \int_0^t e^{-r(s)} \langle\prm \mathbf{A}(s)-\mathcal{A}\bbus, \xms \rangle ds\rightarrow 0,
\end{equation*}
as $m\rightarrow \infty$.
With this we have just shown that the last four terms of \eqref{LIM20} converges to zero as $m\rightarrow \infty$. Then, we can conclude with that
\begin{align}
   \bae e^{-r(\tau_M)} |\xm(\tau_M)|^2+\kappa_1\bae \int_0^{\tau_M}e^{-r(s)} ||\xms||^2_2 ds\rightarrow 0,\label{75*}\\
\bae \int_0^{\tau_M} \int_Z e^{-r(s)} |\prm [\sigma(s,\bu(s),z)-\Sigma(s,z) ]|^2 \eta(dz,ds)\rightarrow 0, \label{LIMXYZ}
%For any $t\in [0,\tau_M]$ $\bu$
\end{align}
as $m\rightarrow \infty$. We easily terminate the proof of the lemma by putting equations \eqref{54-1} and \eqref{LIM17} into \eqref{75*}.
\end{proof}
\begin{proof}[Proof of Proposition \ref{WEAK-COMP2}]
 %\begin{proof}[Proof of Lemma \ref{identity}{\ .}]
First note that for
any $\bw\in V$
\begin{align}
S& =\langle B(\bum,\bum)-B(\bu,\bu),\bw\rangle \nonumber \\
& =\langle B(\bum-\bu,\bum),\bw\rangle +\langle B(\bu,\bum-\bu),\bw\rangle.\label{EQ-B}
\end{align}%
The following equations also hold true
\begin{equation*}
\langle B(\bum-\bu,\bum),\bw\rangle=\langle B(\bum,\bum),\bw\rangle-\langle B%
(\bu,\bum),\bw\rangle,
\end{equation*}%
\begin{equation*}
\langle B(\bum,\bu-\bum),\bw\rangle=\langle B(\bum,\bu),\bw\rangle-\langle B%
(\bum,\bum),\bw\rangle.
\end{equation*}%
Therefore
\begin{equation}
\begin{split}
S=\langle B(\bu,\bu-\bum),\bw\rangle-\langle B(\bum,\bu-\bum),\bw\rangle+\langle B(\bum,\bu),\bw\rangle\\-\langle B(\bu,\bum),\bw\rangle.  \label{CONV18a}
\end{split}
\end{equation}
The operator
\begin{align*}
B_{\mathbf{a},.}:& V\rightarrow V^{\ast } \\
& \bv\mapsto B_{\mathbf{a},.}(\bv)=B(\mathbf{a},\bv)
\end{align*}%
is linear continuous for any fixed $a\in V$. Due to this fact and \eqref{LIM4}, it is true that
\begin{equation}
B(\bu,\bum)\rightharpoonup B(\bu,\bu)\text{ weakly in }%
L^{2}(\Omega\times [0,T];V^{\ast }).  \label{CON16b}
\end{equation}%
By a similar argument, we also prove the following convergence
\begin{equation}
B(\bum,\bu)\rightharpoonup B(\bu,\bu)\text{ weakly in }%
L^{2}(\Omega\times [0,T];V^{\ast }).  \label{CONV19}
\end{equation}%
Now let $\bw$ be an element of $L^{\infty }(\Omega\times[0,T];V)$. We deduce from the property \eqref{B2} that
\begin{equation*}
\begin{split}
& \left\vert \bae \int_{0}^{T}1_{[0,\tau _{M}]}\langle B(\bu(s),\bu(s)-\bums),\bw(s)\rangle-\langle%
B(\bums,\bu(s)-\bums),\bw(s)\rangle ds\right\vert \\
& \leq C\bae \int_{0}^{\tau _{M}}|\bu(s)| ||\bum(s)-\bu(s)||_2
ds+C\bae \int_{0}^{\tau _{M}}|\bum(s)| \,\,||\bum(s)-\bu(s)||_2ds,
\end{split}%
\end{equation*}%
from which and (\ref{LIM13}) we derive that
\begin{equation}
\lim_{m\rightarrow \infty }\bae \int_{0}^{T}1_{[0,\tau _{M}]}\langle B%
(\bu(s),\bu(s)-\bums),\bw(s)\rangle -\langle B(\bums,\bu(s)-\bums),\bw(s)\rangle ds=0  \label{CONV20}
\end{equation}%
Since $\tau _{M}\nearrow T$ almost surely and $L^{\infty }(\Omega\times[0,T];V)$ is dense in $L^{2}(\Omega\times[0,T];V)$, we deduce from (\ref{CONV18a})-(\ref{CONV20})
that the identity \eqref{LIM10} holds.

Next, thanks to the property of $\mathcal{A}_p$ (mainly \eqref{3}) we see that
\begin{equation*}
\bae \int_{0}^{\tau _{M}}\Vert\mathcal{A}_p(\bums)-\mathcal{A}_p(\bu)\Vert_{V^\ast}^{2}ds\leq C \EE \int_{0}^{\tau _{M}}||\bum-\bu||_2^{2}ds.
\end{equation*}%
Owing to (\ref{LIM13}) and the fact that $\tau _{M}\nearrow T$ almost
surely as $M\rightarrow \infty $, we obtain the equation \eqref{LIM11}.

The identity \eqref{LIM12} easily follows from \eqref{LIMXYZ}. This completes the proof of the Proposition \ref{WEAK-COMP2}.
\end{proof}
%\end{proof}
\section{Pathwise uniquness and Convergence of the whole sequence of Galerkin approximation}
In this section we show the pathwise uniqueness of the solution and some (strong) convergences of the Galerkin approximate solution to the
exact solution of \eqref{1}.
\begin{thm}\label{path-Uniq}
 Let $\bu_1$ and $\bu_2$ be two strong solutions to \eqref{2} defined on the same stochastic system $(\Omega, \mathcal{F}, \mathbb{F}, \mathbb{P}, \wie)$. Let $\xi_1$ and $\xi_2$ be their
  respective initial conditions. Then for any $t\in [0,T]$ we have
\begin{equation*}%\label{PU1}
 |\bu_1(t)-\bu_2(t)|^2\le C(\omega) |\xi_1-\xi_2|^2,
\end{equation*}
almost surely.
\end{thm}
\begin{proof}
 Let
 %$(\Omega, \mathcal{F}, \mathbb{F}, \mathbb{P}, \wie)$ be a stochastic basis (!! \textbf{Definition}) and
 $\bu_1$ (resp., $\bu_2$) be a strong solution to \eqref{2}
 with initial condition $\xi_1$ (resp., $\xi_2$). Let $\bw=\bu_1-\bu_2$ and $\xi=\xi_1-\xi_2$. It is not hard to see that
\begin{equation*}
\begin{split}
 \bw(t)+\kappa_1\int_0^t A\bw(s) ds+\kappa_0 \int_0^t \left(\mathcal{A}_p \bu_1(s) -\mathcal{A}_p \bu_2(s) \right)ds\\
=\xi+\int_0^t \int_Z \left( \sigma(s,\bu_1(s),z)-\sigma(s,\bu_2(s),z) \right) \wie(dz,ds)\\
+\int_0^t \left(B(\bu_2(s), \bu_2(s))-B(\bu_1(s), \bu_1(s))\right) ds.
\end{split}
\end{equation*}
Applying It\^o's formula to the function $\Phi(x)=|x|^2$ and $\bw(t)$ implies that
\begin{equation*}
 \begin{split}
  |\bw(t)|^2+2 \kappa_1\int_0^t ||\bw(s)||^2_2 ds+2 \kappa_0 \int_0^t \langle \mathcal{A}_p \bu_1(s) -\mathcal{A}_p \bu_2(s), \bw(s) \rangle ds\\
=|\xi|^2+2 \int_0^t \int_Z \left( \sigma(s,\bu_1(s),z)-\sigma(s,\bu_2(s),z), \bw(s) \right) \wie(dz,ds)\\
-2 \int_0^t \langle B(\bu_1(s), \bu_1(s))-B(\bu_2(s), \bu_2(s)), \bw(s)\rangle ds\\
+\int_0^t \int_Z |\sigma(s,\bu_1(s),z)-\sigma(s,\bu_2(s),z)|^2 \eta(dz,ds).
 \end{split}
\end{equation*}
Next we introduce the real valued process $$\rho(t)=e^{-\frac{C_0^2}{\kappa_1}\int_0^t ||\bu_1(s)||^2_2 ds}.$$ Now we apply It\^o's formula to $\rho(t)|\bw(t)|^2$ and we get
\begin{equation*}
 \begin{split}
   \rho(t) |\bw(t)|^2+2 \kappa_1\int_0^t \rho(s)||\bw(s)||^2_2 ds+2 \kappa_0 \int_0^t \rho(s) \langle \mathcal{A}_p \bu_1(s) -\mathcal{A}_p \bu_2(s), \bw(s) \rangle ds\\
=-2 \int_0^t \rho(s)\langle B(\bu_1(s), \bu_1(s))-B(\bu_2(s), \bu_2(s)), \bw(s)\rangle ds-\frac{C_0^2}{\kappa_1}\int_0^t \rho(s) |\bw(s)|^2 ||\bu_1(s)||^2_2 ds \\
+2 \int_0^t \int_Z \rho(s)\left( \sigma(s,\bu_1(s),z)-\sigma(s,\bu_2(s),z), \bw(s) \right) \wie(dz,ds)\\
+\int_0^t \int_Z \rho(s) |\sigma(s,\bu_1(s),z)-\sigma(s,\bu_2(s),z)|^2 \eta(dz,ds)+|\xi|^2.
 \end{split}
\end{equation*}

By making use of \eqref{AP3}, \eqref{EQ-B}, \eqref{B1}, \eqref{B2} and Young's inequality with $\eps=\kappa_1$ in the above estimate and by taking the mathematical expectation to both
sides of the resulting estimate yield
\begin{equation*}
 \begin{split}
  \bae \rho(t)|\bw(t)|^2+2 \kappa_1\bae \int_0^t\rho(s) ||\bw(s)||^2_2 ds\le \bae  \int_0^t \int_Z \rho(s) |\sigma(s,\bu_1(s),z)-\sigma(s,\bu_2(s),z)|^2 \nu(dz)ds\\
+|\xi|^2
.
 \end{split}
\end{equation*}
Using point (1) of Condition \ref{cond-s-1} yields that
\begin{equation*}
 \bae \rho(t)|\bw(t)|^2 \le \bae |\xi|^2 + \bae \int_0^t \rho(s) |\bw(s)|^2 ds,
\end{equation*}
from which and Gronwall's Lemma we deduce the existence ofa constant $C>0$ such that
\begin{equation*}
 \bae \rho(t)|\bw(t)|^2 \le C\bae |\xi|^2,
\end{equation*}
for any $t\in [0,T]$.
Since $\rho(t)$ is bounded $\PP-$a.s, we conclude easily from the last estimate the proof of the theorem.

\end{proof}
Next we will show that the whole sequence of solutions to the Galerkin approximation system \eqref{discrete-sys}
converges in mean square to the exact strong solution of \eqref{1}. Mainly we have
\begin{thm}\label{STRONGCONV}
 The whole sequence of Galerkin approximation $\{\bum: m\in \NN\}$ defined by \eqref{discrete-sys} satisfies
\begin{align}
 \lim_{m\rightarrow \infty}\bae |\bum(T-)-\bu(T-)|^2=0,\label{CV1}\\
\lim_{m\rightarrow \infty}\bae\int_0^{T-} \Vert \bums-\bu(s)\Vert^2_2ds=0.\label{CV2}
\end{align}
\end{thm}
%The main ingredient of the proof of this result is the following lemma whose proof can be found in \cite{Breckner}.
The main ingredient of the proof of this result is the following lemma, its proof follows a very small modification of the proof of
\cite[Proposition B.3]{breckner}.
\begin{lem}
 Let $\{Q_m; m\ge 1\}\subset L^2(\Omega\times [0,T];\mathbb{R}))$  be a sequence of \cadlag real-valued process, and let
$\{T_M;M\ge 1\}$ be a sequence of $\mathcal{F}^t$-stopping times such that $T_M$ is increasing to $T$, $\sup_{m\ge 1}\mathbb{E}|Q_m(T)|^2<\infty$, and
$\lim_{m\rightarrow \infty}\mathbb{E}|Q_m(T_M)|=0$ for all $M\ge 1$. Then $\lim_{m\rightarrow \infty}\mathbb{E}|Q_m(T-)|=0$.
\end{lem}
\begin{proof}[ Proof of Theorem \ref{STRONGCONV}]

It follows from Lemma \ref{COMP} that
\begin{equation}\label{STRONGCONV1}
 \lim_{m\rightarrow \infty}\mathbb{E}\int_0^{\tau_M}||\bum(t)-\bu(t)||^2_2 dt=0,
\end{equation}
and
\begin{equation}\label{STRONGCONV2}
 \lim_{m\rightarrow \infty}\mathbb{E}|\bum(\tau_M)-\bu(\tau_M)|^2=0,
\end{equation}
for any $M\ge 1$.
So by applying the preceding lemma to $Q_m(t)=|\bum(t)-\bu(t)|^2$, $T_M=\tau_M$ and taking into account
\eqref{STRONGCONV2},
the estimates in Proposition \ref{proposition-estimate} and the uniqueness of $\bu$,  we see that the whole sequence $\bum$ defined by
\eqref{discrete-sys} satisfies \eqref{CV1}. To prove \eqref{CV2} we need an extra estimate for the sequence  $\{\bum: m \in \mathbb{N}\}$. Since
$$\int_0^t \int_Z \Psi(s,z)\tilde{\eta}(dz,ds)=\int_0^t \int_Z \Psi(s,z)\eta(ds,dz)-\int_0^t \int_Z \Psi(s,z)\nu(dz)ds,$$ and
$|\bv|^2+2 (\bv,\bw)=|\bv+\bw|^2-|\bw|^2$ we deduce from \eqref{eq17} that
\begin{equation*}
\begin{split}
  4 \kappa_1^2\EE \biggl(\int_0^t \lVert \bums\rVert^2_2 ds\biggr)^2\le& 2 \EE \lvert \xi\rvert^4+
 4 \EE \biggl(\int_0^t \int_Z \lvert \sigma(s, \bums, z)\rvert \nu(dz)ds\biggr)^2
 \\ &\quad + 4 \EE\biggl(\int_0^t\int_Z \biggl[\lvert \sigma(s,\bums,z)+\buns\rvert^2-\lvert \buns\rvert^2 \biggr]\tilde{\eta}(dz,ds)\biggr)^2.
\end{split}
\end{equation*}
Using item (1) of Condition \ref{cond-s-1} and the estimate in Proposition \ref{proposition-estimate} we infer from the last inequality that
\begin{equation*}
 \begin{split}
  4 \kappa_1^2\EE \biggl(\int_0^t \lVert \bums\rVert^2_2 ds\biggr)^2\le 4 \EE\biggl(\int_0^t\int_Z \biggl[\lvert \sigma(s,\bums,z)+\buns\rvert^2-\lvert \buns\rvert^2 \biggr]\tilde{\eta}(dz,ds)\biggr)^2
  \\+2 \EE \lvert \xi\rvert^4+ 4 \ell_0^2 T^2+4C \ell_1^2 T^2 (\EE \lvert \xi\rvert^4 +1).
 \end{split}
\end{equation*}
Now invoking \cite[Theorem 4.14]{Rudiger} we see that
\begin{equation*}
 \begin{split}
  4 \kappa_1^2\EE \biggl(\int_0^t \lVert \bums\rVert^2_2 ds\biggr)^2\le 4 \EE\int_0^t\int_Z \biggl[\biggl\lvert \lvert \sigma(s,\bums,z)+\buns\rvert^2-\lvert \buns\rvert^2 \biggr\rvert^2\biggr]\nu(dz)ds
  \\+2 \EE \lvert \xi\rvert^4+ 4 \ell_0^2 T^2+4C \ell_1^2 T^2 (\EE \lvert \xi\rvert^4 +1),
 \end{split}
\end{equation*}
from which with item (3) of Condition \ref{cond-s-1} and Proposition \ref{proposition-estimate} we derive that
\begin{equation*}
 4 \kappa_1^2\EE \biggl(\int_0^t \lVert \bums\rVert^2_2 ds\biggr)^2\le \left(2+4C\ell_1^2 T^2+C\ell_3T+CT \right)\EE \lvert \xi\rvert^4
 + 4\ell_0^2 T^2+4C \ell_1^2 T^2+CT.
\end{equation*}
This also implies that
\begin{equation*}
 4 \kappa_1^2\EE \biggl(\int_0^t \lVert \bu(s)\rVert^2_2 ds\biggr)^2\le \left(2+4C\ell_1^2 T^2+C\ell_3T+CT \right)\EE \lvert \xi\rvert^4
 + 4\ell_0^2 T^2+4C \ell_1^2 T^2+CT.
\end{equation*}
We see easily from the last two estimates and \eqref{STRONGCONV1} that $Q_m(t)=\int_0^t||\bum(s)-\bu(s)||_2^2 ds$, $T_M=\tau_M$ satisfy
the hypotheses
of the above lemma, therefore we can deduce that \eqref{CV2} holds.
This ends the proof of Theorem \ref{STRONGCONV}.
\end{proof}
\section{Existence and ergodicity of invariant measure}

In this section we are intereted in the study of some qualitative properties of the solution of \eqref{1}. We will mainly analyse
the Markov, Fellerian properties of the solution. We will also derive the existence of ergodic invariant measures.
To start with our investigation we denote  by $\bu(t;\xi)$ the solution of \eqref{1} with initial
condition $\xi\in H$, and  by $C_b(H)$ we describe the space of all continuous real-valued functionals defined on $H$. Next, we define the family of mappings $\{\CP_t, t\ge 0\}$ ($\CP_t$ for short ) defined on $C_b(H)$ by
\begin{equation*}
 \CP_t \phi(\xi)=\EE\phi(\bu(t;\xi)),
\end{equation*}
for any $\phi\in C_b(H)$, $\xi\in H$, and $t\ge 0$. Thanks to Theorem \ref{main-ex} and Theorem \ref{path-Uniq} the family $\CP_t$ defines a semigroup on
$C_b(H)$.
More properties of the solution $\bu(t;\xi)$ and the semigroup $\CP_t$ are given in the following results.
\begin{thm}\label{Feller}
 The solution $\bu(t;\xi)$ defines a Markov process
%  with probability transition $\CP_t$ defined by
%  $$ \CP_t(A)=\EE \mathbf{1}_A(\bu(t;\xi)),$$
%  for any $A\in \mathcal{B}(H)$.
%
and the semigroup $\CP_t$ is Fellerian; that is,  $\CP_t$ satisfies the following
$$ \CP_t \left(C_b(H)\right)\subset C_b(H),$$ for any $t\ge 0$.
\end{thm}
Before we proceed to the proof of these statements let us give an auxillary result.
\begin{lem}\label{Feller-1}
Let $\bu(t;\xi_1)$ and $\bu(t;\xi_2)$ be two solutions of \eqref{1} associated to two distincts initial conditions $\xi_1$, $\xi_2$, let
 \begin{equation}\label{taur}
  \tau_R^\xi=\inf\{t: |\bu(t;\xi)|> R\}, \forall R>0, \xi\in H.
 \end{equation}
 Let us set $\tau^{\xi_1,\xi_2}_R=\tau_R^{\xi_1}\wedge \tau_R^{\xi_2}$, $t_R=t\wedge \tau_R^{\xi_1,\xi_2}$ and $\bw(t)=\bu(t;\xi_1)-\bu(t;\xi_2)$, $t\in [0,\infty)$. Then for any $R>0$ and
  $t\in [0,\infty)$ there exists a positive constant $C$ such that
 \begin{equation}\label{eqFeller-1}
  \EE \lvert \bw(t_R)\rvert^2\le C \EE \lvert \xi_1-\xi_2\rvert^2.
 \end{equation}
\end{lem}
\begin{proof}[Proof of Lemma \ref{Feller-1}]
 As in the proof of Theorem \ref{path-Uniq} we can check by making use of It\^o's formula that $\lvert\bw(t_R)\rvert^2$ satisfies
 \begin{equation*}
 \begin{split}
  \lvert\bw(t_R)\rvert^2+2\kappa_1\int_0^{t_R} \lVert \bu(s;\xi)\rVert^2_2 ds\le \lvert \xi_1-\xi_2\rvert^2+2 \int_0^{t_R} \langle B(\bw(s), \bu(s;\xi_1), \bw(s)\rangle ds
 \\ +\int_0^{t_R} \int_Z |\sigma(s, \bu(s;\xi_1),z)-\sigma(s, \bu(s;\xi_2),z)|^2\eta(dz,ds)\\+
  2\int_0^{t_R}\int_Z \left(\sigma(s, \bu(s;\xi_1),z)-\sigma(s, \bu(s;\xi_2),z),\bw(s-)\right) \tilde{\eta}(dz,ds).
 \end{split}
 \end{equation*}
 Using the skew-symmetricity of $B$ and H\"older's inequality we derive from the last inequality that
  \begin{equation*}
 \begin{split}
  \lvert\bw(t_R)\rvert^2+2\kappa_1\int_0^{t_R} \lVert \bu(s;\xi)\rVert^2_2 ds\le \lvert \xi_1-\xi_2\rvert^2
  +2 C \int_0^{t_R}\biggl(\lvert \bw(s)\cdot \nabla \bw(s)\rvert \times \lvert\bu(s;\xi_1 )\rvert \biggr) ds
 \\ +\int_0^{t_R} \int_Z |\sigma(s, \bu(s;\xi_1),z)-\sigma(s, \bu(s;\xi_2),z)|^2\eta(dz,ds)\\+
  2\int_0^{t_R}\int_Z \left(\sigma(s, \bu(s;\xi_1),z)-\sigma(s, \bu(s;\xi_2),z),\bw(s-)\right) \tilde{\eta}(dz,ds).
 \end{split}
 \end{equation*}
 Owing to H\"older's inequality and the fact $\lvert \bu(s;\xi_1)\rvert\ge R$ on $[0,t_R]$ we infer the existence of a constant $C_R=C(R)>0$ such that
\begin{equation*}
 \begin{split}
  \lvert\bw(t_R)\rvert^2+2\kappa_1\int_0^{t_R} \lVert \bu(s;\xi)\rVert^2_2 ds\le \lvert \xi_1-\xi_2\rvert^2
  +2 C_R \int_0^{t_R}\biggl(\lvert \bw(s)\rvert \times \lvert \nabla \bw(s)\rvert_{L^q} \biggr) ds
 \\ +\int_0^{t_R} \int_Z |\sigma(s, \bu(s;\xi_1),z)-\sigma(s, \bu(s;\xi_2),z)|^2\eta(dz,ds)\\+
  2\int_0^{t_R}\int_Z \left(\sigma(s, \bu(s;\xi_1),z)-\sigma(s, \bu(s;\xi_2),z),\bw(s-)\right) \tilde{\eta}(dz,ds),
 \end{split}
 \end{equation*}
 where $2<q\le \frac{2n}{n-2}$.
 Thanks to Young's inequality and the continuous embedding $\mathcal{H}^1\subset L^q$ we easily see that
 \begin{equation*}
 \begin{split}
  \lvert\bw(t_R)\rvert^2+2\kappa_1\int_0^{t_R} \lVert \bu(s;\xi)\rVert^2_2 ds\le \lvert \xi_1-\xi_2\rvert^2
  +\frac{2 C_R}{\eps} \int_0^{t_R}\lvert \bw(s)\rvert^2 ds+\eps \int_0^{t_R}  \lVert \bw(s)\rVert_2^2ds
 \\ +\int_0^{t_R} \int_Z |\sigma(s, \bu(s;\xi_1),z)-\sigma(s, \bu(s;\xi_2),z)|^2\eta(dz,ds)\\+
  2\int_0^{t_R}\int_Z \left(\sigma(s, \bu(s;\xi_1),z)-\sigma(s, \bu(s;\xi_2),z),\bw(s-)\right) \tilde{\eta}(dz,ds),
 \end{split}
 \end{equation*}
 Choosing $\eps=\kappa_1$, using item (2) of Condition \ref{cond-s-1} and taking the mathematical expectation yield that
  \begin{equation}\label{Feller-2}
 \begin{split}
 \EE \lvert\bw(t_R)\rvert^2+\kappa_1\EE \int_0^{t_R} \lVert \bu(s;\xi)\rVert^2_2 ds\le \lvert \xi_1-\xi_2\rvert^2
  +\left(\frac{2 C_R}{\kappa_1}+ L_1\right) \int_0^{t_R}\lvert \bw(s)\rvert^2 ds.
 \end{split}
 \end{equation}
 where we have used the fact that
 \begin{align*}
  \EE \int_0^{t_R} \int_Z |\sigma(s, \bu(s;\xi_1),z)-\sigma(s, \bu(s;\xi_2),z)|^2\eta(dz,ds)\\
  =\EE \int_0^{t_R} |\sigma(s, \bu(s;\xi_1),z)-\sigma(s, \bu(s;\xi_2),z)|^2\nu(dz)ds,\\
  \end{align*}
  and
  \begin{align*}
  2\EE\int_0^{t_R}\int_Z \left(\sigma(s, \bu(s;\xi_1),z)-\sigma(s, \bu(s;\xi_2),z),\bw(s-)\right) \tilde{\eta}(dz,ds)=0.
 \end{align*}
Notice that \eqref{Feller-2} can be rewritten in the following form
\begin{equation*}%\label{Feller-3}
 \begin{split}
 \EE \lvert\bw(t_R)\rvert^2+\kappa_1\EE \int_0^{t_R} \lVert \bu(s;\xi)\rVert^2_2 ds\le \lvert \xi_1-\xi_2\rvert^2
  +\left(\frac{2 C_R}{\kappa_1}+ L_1\right) \int_0^{t}\lvert \bw(s\wedge \tau_R)\rvert^2 ds,
 \end{split}
 \end{equation*}
 from which along with the application Gronwall's lemma we deduce the existence of a positive constant $C=C(t,R)$ such that
 \begin{equation*}
  \EE \lvert\bw(t_R)\rvert^2\le C \lvert \xi_1-\xi_2\rvert^2.
 \end{equation*}
The proof of the lemma is now finished.
\end{proof}
Now we continue with the proof of Theorem \ref{Feller}
\begin{proof}[Proof of Theorem \ref{Feller}]
 Owing to the Theorem \ref{path-Uniq} and the fact that $\tilde{\eta}(A\times [0,t])$, $A\times [0,t]\in \mathcal{B}(Z\times \mathbb{R}_{+})$ is time homogeneous,
 the Markovian property of $\bu(t;\xi)$, $\xi\in H$, can be checked using the same argument as in \cite{daprato} (see also
 \cite{Albeverio+Mandrekar+Rudiger}). Now we want to check that
 $\CP_t\left(C_b(H)\right)\subset C_b(H) $. For this purpose let us consider $\xi\in H$ and a sequence $\{\xi_m: m\in \mathbb{N}\}\subset H$ such that
 $\xi_m\rightarrow \xi$ as $m\rightarrow \infty$. Let us prove that
 \begin{equation*}
  \CP_t\phi(\xi_m)\rightarrow \CP_t \phi(\xi), \forall \phi\in C_b(H),
 \end{equation*}
as $m$ tends to infinity. To shorten notation we set $\tau_R=\tau_R^{\xi_m}\wedge \tau_R^{\xi}$ where the stopping times $\tau_R^{\xi}$ is defined as in \eqref{taur}.
For any $t\in [0,T], T\ge 0$ and $\phi\in C_b(H)$, we have
\begin{equation*}
 \begin{split}
  \lvert \CP_t\phi(\xi_m)-\CP_t \phi(\xi)&=\biggl \lvert \EE\biggl(\biggl[\phi(\bu(t;\xi_m))-\phi(\bu(t;\xi))\biggr]\mathbf{1}_{[t<\tau_R]\cup [t\ge \tau_R]}\biggr)\biggr\rvert,\\
   &\le \biggl \lvert \EE\biggl(\biggl[\phi(\bu(t;\xi_m))-\phi(\bu(t;\xi))\biggr]\left(\mathbf{1}_{[t\ge \tau^{\xi_m}_R]}+\mathbf{1}_{[t\ge \tau^{\xi}_R]}\right)\biggl)\biggr\rvert
   \\
   & \quad +\biggl \lvert \EE\biggl(\biggl[\phi(\bu(t;\xi_m))-\phi(\bu(t;\xi))\biggr]\mathbf{1}_{[t<\tau_R]}\biggr)\biggr \rvert.
 \end{split}
\end{equation*}
Thanks to the fact that $\EE \lvert \bu(t;\xi)\rvert^2< C(\xi), \forall \xi\in H$ (see the estimate in Theorem \ref{main-ex}), we obtain that for any $\eps>0$ there exists $m_1$ such
that for any $R> m_1$
\begin{equation*}
 \PP\left(\tau_R^{\xi_m}\ge t\right)+\PP \left(\tau_R^{\xi}\ge t\right)\le \frac{\eps}{4 \lVert \phi\rVert_{\infty}},
\end{equation*}
where $$\lVert \phi\rVert_\infty=\sup_{\xi\in H} |\phi(x)|.$$
Thus
\begin{equation*}
%\begin{multline}
 \lvert \CP_t\phi(\xi_m)-\CP_t\phi(\xi)\rvert\le \biggl \lvert \EE\biggl(\biggl[\phi(\bu(t;\xi_m))-\phi(\bu(t;\xi))\biggr]\mathbf{1}_{[t<\tau_R]}\biggr)\biggr \rvert+
 2 \lVert \phi\rVert_\infty \frac{\eps}{4\lVert \phi\rVert_\infty}.
% \end{multline}
\end{equation*}
That is,
\begin{equation*}
%\begin{multline}
 \lvert \CP_t\phi(\xi_m)-\CP_t\phi(\xi)\rvert\le \biggl \lvert \EE\biggl(\biggl[\phi(\bu(t;\xi_m))-\phi(\bu(t;\xi))\biggr]\mathbf{1}_{[t<\tau_R]}\biggr)\biggr \rvert+
 \frac{\eps}{2}.
% \end{multline}
\end{equation*}
Since $\mathbf{1}_{[t<\tau_R]}\le 1$ and $t\wedge \tau_R=t$ when $t<\tau_R$, we readily have that
\begin{equation}\label{Feller-4}
  \lvert \CP_t\phi(\xi_m)-\CP_t\phi(\xi)\rvert\le \biggl \lvert \EE\biggl(\biggl[\phi(\bu(t_R;\xi_m))-\phi(\bu(t_R;\xi))\biggr]\biggr)\biggr \rvert+
 \frac{\eps}{2},
\end{equation}
where we have put $t_R=t\wedge \tau_R$. By the continuuity of $\phi$, for the same $\eps>0$ as above we can find $\kappa>0$ such that if
$\lvert \bu(t_R; \xi_m)-\bu(t_R; \xi)\rvert< \kappa$ we have
\begin{equation}\label{F2}
 \lvert \phi(\bu(t_R;\xi_m))-\phi(\bu(t_R;\xi))\rvert< \frac{\eps}{4}.
\end{equation}
Note that from \eqref{Feller-4} we derive that
\begin{equation*}
 \begin{split}
  \lvert \CP_t\phi(\xi_m)-\CP_t\phi(\xi)\rvert\le \biggl \lvert
  \EE\biggl(\biggl[\phi(\bu(t_R;\xi_m))-\phi(\bu(t_R;\xi))\biggr]\mathbf{1}_{\{\lvert\bu(t_R; \xi_m)-\bu(t_R; \xi) \rvert\ge \kappa\}}\biggr)\biggr \rvert\\
 +\EE\biggl(\biggl[\phi(\bu(t_R;\xi_m))-\phi(\bu(t_R;\xi))\biggr]\mathbf{1}_{\{\lvert\bu(t_R; \xi_m)-\bu(t_R; \xi) \rvert<\kappa\}}\biggr)\biggr \rvert+\frac{\eps}{2},
 \end{split}
\end{equation*}
from which all together with \eqref{F2} we derive that
\begin{equation}\label{F1}
 \begin{split}
  \lvert \CP_t\phi(\xi_m)-\CP_t\phi(\xi)\rvert\le 2 \lVert \phi\rVert_\infty\PP\biggl(\lvert\bu(t_R; \xi_m)-\bu(t_R; \xi) \rvert\ge \kappa\biggr)\\
  ++\EE\biggl(\biggl[\phi(\bu(t_R;\xi_m))-\phi(\bu(t_R;\xi))\biggr]\mathbf{1}_{\{\lvert\bu(t_R; \xi_m)-\bu(t_R; \xi) \rvert<\kappa\}}\biggr)\biggr \rvert+\frac{\eps}{2}.
 \end{split}
\end{equation}
Invoking the estimate \eqref{eqFeller-1} and Chebychev's ineqlity we obtain that
\begin{equation}\label{F1-b}
  2 \lVert \phi\rVert_\infty\PP\biggl(\lvert\bu(t_R; \xi_m)-\bu(t_R; \xi) \rvert\ge \kappa\biggr)\le \frac{2 \lVert \phi\rVert_\infty C}{\kappa^2}\lvert \xi_m-\xi\rvert^2.
\end{equation}
But as $\xi_m\rightarrow \xi$ as $m\rightarrow \infty$ we have that for any $\delta>0$ there exists $m_2>0$ such that if $m>m_2$ we have $\lvert\xi_m-\xi\rvert^2<\delta$. Choosing
$\delta=\frac{\eps \kappa^2}{8C\lVert \phi\rVert_\infty}$ we can derive from \eqref{F1-b} that
\begin{equation}\label{F3}
  2 \lVert \phi\rVert_\infty\PP\biggl(\lvert\bu(t_R; \xi_m)-\bu(t_R; \xi) \rvert\ge \kappa\biggr)\le \frac{\eps}{4}.
\end{equation}
So combining \eqref{F2}, \eqref{F1} and \eqref{F3} we see that for any $\eps>0$ there exists $m_0>0$ such that if $m> m_0$ then
\begin{equation*}
 \lvert \CP_t\phi(\xi_m)-\CP_t\phi(\xi)\rvert< \eps,
\end{equation*}
which shows that $\CP_t$ is a Fellerian semigroup.
\end{proof}
Owing to Theorem \ref{Feller} we can discuss about the existence of the invariant measure associated to the semigroup $\CP_t$.
\begin{thm}\label{mes-inv}
 The Markovian semigroup $\CP_t$ has at least one invariant measure $\mu$. Moreover, $\mu$ is concetrated on $V$, i.e, $\mu(V)=1$.
\end{thm}
\begin{proof}
 Let $\{T_n; n\in \mathbb{N}\}\subset [0,\infty)$ be a sequence such that $T_n\nearrow \infty$ as $n\rightarrow \infty$. For any $A\in \mathcal{B}(H)$ let us set
 \begin{equation*}
  \mu_n(A)=\frac{1}{T_n}\int_0^{T_n} \PP\left(\bu(t;\xi)\in A\right) dt.
 \end{equation*}
It is clear that $\mu_n$ defines a measure on $(H, \mathcal{B}(H))$. Let $R>0$ and $A_R=\{\bu: \lVert \bu\rVert_2 > R\}$. Using Chebychev's inequality anf Fubini's Theorem
we see that
\begin{equation*}
 \mu_n(A_R)\le \frac{1}{R^2}\frac{1}{T_n}\EE\int_0^{T_n} \lVert \bu(s;\xi) \rVert^2_2ds.
\end{equation*}
Owing to the estimate in Theorem \ref{main-ex} we have that
\begin{equation*}
 \mu_n(A_R)\le \frac{C(1+\lvert \xi\rvert^2)}{R^2}.
\end{equation*}
This implies that $\mu_n(A_R)\rightarrow 0$ uniformly in $n$ as $R\rightarrow \infty$. Since the ball $B_R=V\backslash A_R$ is compact in $H$, we conclude that the family of measures
$\mu_n$ is tight on $H$. This yields that there exists a subsequence $\mu_{n_k}$ and a measure $\mu$ defined on $(H, \mathcal{B}(H))$ such that
\begin{equation*}
 \int_H \phi(x) \mu_{n_k}(dx)\rightarrow \int_H \phi(x) \mu(dx), \forall \phi\in C_b(H).
\end{equation*}
Since $\CP_t$ satisfies the Markov-Feller property, we can infer from Krylov-Bogoluibov's theorem that it admits an invariant measure which is equal to $\mu$.

It remains to show that $\mu$ is concentrated on $V$. For this purpose it is sufficient to show that $\mu(H\backslash V)=0$.
To do so we will first show that
$$ \mu_n(H\backslash V)=0, \forall n.$$ Thanks to the estimate in Theorem \ref{main-ex} we can find a set $I\times \Omega_0\subset \Omega_{T_n}, T_n\ge 0$ ($\Omega_{T_n}=[0, T_n]\times \Omega$) with
$\lambda\otimes \PP(\Omega_t\backslash I\times \Omega_0)=0$ and $\bu(t;\xi)(\omega)\in V$ for any $(t, \omega)\in I\times \Omega_0$. This fact implies that
\begin{equation*}
 \PP\biggl(\int_0^{T_n} \mathbf{1}_N (t,\omega) dt\biggr)=0,
\end{equation*}
where $$ N=\{(t,\omega)\in \Omega_{T_n}: \bu(t;\xi)(\omega)\in H\backslash V\}.$$
Owing to Fubini's theorem we infer the existence of $J\subset [0, T_n]$ with $\lambda([0, T_n]\backslash J)=0$ and
\begin{equation*}
 \PP\left(\{\omega \in \Omega: \bu(t;\xi) \in H\backslash V\}\right)=0,
\end{equation*}
for any $t\in J$. Setting $N_t=\{\omega\in \Omega; \bu(t;\xi) \in H\backslash V\}$ for any $t\in J$, we find that
\begin{align*}
 \mu_n(H\backslash V)&=\frac{1}{T_n}\int_0^{T_n} \PP(N_t)dt,\\
 &=\frac{1}{T_n}\int_0^{T_n} \mathbf{1}_{J}(t) \PP(N_t) dt,\\
 &=0.
\end{align*}
This means that the support of $\mu_n$ is included in $V$. Since $\mu$ is the weak limit of $\mu_n$, we derive from {\color{red} \bf \cite[Theorem 2.2]{Chow+Kashminskii}} that the support of
$\mu$ is included in $V$.
\end{proof}
Our next concern is to check whether the invariant measure $\mu$ is ergodic or not. In fact we will find that it is ergodic provided that $\kappa_1$ is large enough.
We will make our claim clearer later on, but for now let us prove an important fact about the invariant measure $\mu$.
\begin{prop}
 If $2\kappa_1\lambda_1^2-\ell_1>0$, then  there exists a constant $\tilde{L}>0$ depending only on $\kappa_1, \lambda_1, \ell_0, \ell_1$ such that
 \begin{equation}
  \int_H \left(\lvert \xi\rvert^2 +\lVert \xi\rVert_2^2 \right)\mu(dx) < \tilde{L}.
 \end{equation}
\end{prop}
\begin{proof}
 First we should notice that by It\^o's formula we have
 \begin{equation}\label{inv1}
 \begin{split}
  & \lvert \bu(t;\xi)\rvert^2+2\kappa_1 \int_0^t \lVert \bu(s;\xi)\rVert^2_2 ds+ 2\int_0^t \langle \mathcal{A}_p \bu(s;\xi), \bu(s,\xi)\rangle ds\\
  &\quad\quad =\lvert \xi\rvert^2 + \int_0^t \int_Z \lvert \sigma(s, \bu(s; \xi), z)\rvert^2 \eta(dz,ds)
  + 2\int_0^t \int_Z \left(\sigma(s, \bu(s;\xi),z) ,\bu(s-;\xi) ) \right)\tilde{\eta}(dz,ds).
 \end{split}
 \end{equation}
Now for any $\eps>0$ let $\Phi(y)=\frac{y}{1+\eps y }, y\in \mathbb{R}_{+}$. It is clear that
\begin{align*}
 \Phi^\prime (y)&=\frac{1}{(1+\eps y)^2},\\
 \Phi^{\prime\prime}(y)&=\frac{-2\eps }{(1+\eps y)^3},
\end{align*}
for any $y\ge 0$. It is clear from the last equality that $\Phi^{\prime \prime}(y)< 0$, and  $|\Phi^{\prime \prime}|\le 2 \eps$ for any $y\ge 0.$
Notice also that $\eta(dz,ds)=\tilde{\eta}(dz,ds)+\nu(dz)ds$  and
$$\lvert \sigma(s, \bu(s,\xi), z)\rvert^2+2\left(\sigma(s, \bu(s,\xi), z), \bu(s-;\xi) \right)=\lvert \sigma(s, \bu(s,\xi), z)+\bu(s-;\xi)\rvert^2-\lvert\bu(s-;\xi)\rvert^2. $$
By setting $Y(t)=\lvert \bu(t;\xi)\rvert^2 $ and $\Psi=\lvert \sigma(s, \bu(s,\xi), z)+\bu(s-;\xi)\rvert^2-\lvert\bu(s-;\xi)\rvert^2$ we can rewrite \eqref{inv1} in the following form
\begin{equation*}
\begin{split}
 Y(t)+2\kappa_1 \int_0^t \lVert \bu(s;\xi)\rVert^2_2 ds+2 \int_0^t \langle \mathcal{A}_p\bu(s;\xi), \bu(s;\xi)\rangle ds=\lvert\xi\rvert^2
 + \int_0^t \int_Z \lvert \sigma(s, \bu(s;\xi), z) \rvert^2 \nu(dz)ds \\
 +\int_0^t \int_Z \Psi \tilde{\eta}(dz,ds).
\end{split}
\end{equation*}
Applying It\^o's formula to $\Phi(Y)$ we obtain that
\begin{equation*}
 \begin{split}
  \Phi(Y(t) )+2\kappa_1 \int_0^t \Phi^\prime (Y(s))\lVert \bu(s;\xi)\rVert^2_2 ds+2 \int_0^t \Phi^\prime (Y(s)) \langle \mathcal{A}_p\bu(s;\xi), \bu(s;\xi)\rangle ds
 \\ =\Phi(\lvert\xi\rvert^2)  +\int_0^t\int_Z \biggl(\Phi(Y(s-)+\Psi)-\Phi(Y(s-))-\Phi^\prime(Y(s-))\Psi \biggr)\eta(dz,ds)\\
  + \int_0^t \Phi^\prime (Y(s)) \int_Z \lvert \sigma(s, \bu(s;\xi), z) \rvert^2 \nu(dz)ds
 + \int_0^t \int_Z \biggl(\Phi(Y(s-) +\Psi)-\Phi(Y(s-))\biggr)\tilde{\eta}(dz,ds).
 \end{split}
\end{equation*}
Since $\langle \mathcal{A}_p\bu(s;\xi), \bu(s;\xi)\rangle\ge 0$ and $\Phi^\prime(y)>0$ for any $y\ge 0$, we can drop out the third term from the left-hand side of the last equation.
Therefore
we obtain that
\begin{equation*}%\label{inv2}
 \begin{split}
   \Phi(Y(t) )+2\kappa_1 \int_0^t \Phi^\prime (Y(s))\lVert \bu(s;\xi)\rVert^2_2 ds
  \le \Phi(\lvert\xi\rvert^2)  + \int_0^t \Phi^\prime (Y(s)) \biggl(\int_Z \lvert \sigma(s, \bu(s;\xi), z) \rvert^2 \nu(dz)\biggr)ds\\
  +\int_0^t\int_Z \biggl(\int_0^1 \Phi^{\prime \prime}(Y(s-)+\theta \Psi)\Psi^2 d\theta \biggr)\left(\tilde{\eta}(dz,ds)+\nu(dz)ds\right)\\
 + \int_0^t \int_Z \biggl(\int_0^1 \Phi^\prime(Y(s-)+\theta \Psi)\Psi d\theta\biggr)\tilde{\eta}(dz,ds),
 \end{split}
\end{equation*}
where we have used the identities
\begin{align*}
 \Phi(y+\psi)-\Phi(y)=\int_0^1 \Phi^\prime(y+\theta \psi)\psi d\theta,\\
 \Phi(y+\psi)-\Phi(y)-\Phi^\prime(y)\psi=\int_0^1 \Phi^{\prime\prime}(y+\theta \psi)\psi^2 d\theta.
\end{align*}
Since $\lvert \Phi^\prime(\cdot)\rvert<1$ and $\lvert \Phi^{\prime\prime}\rvert<2\eps$ and
\begin{equation*}
 \EE \Psi^{r}\le C  \EE (1+\lvert\bu(s;\xi)\rvert^{2r})<C,
\end{equation*}
with $r=1,2$, the stochastic integrals
\begin{align*}
 \int_0^t\int_Z \biggl(\int_0^1 \Phi^{\prime}(Y(s)+\theta \Psi)\Psi d\theta \biggr)\tilde{\eta}(dz,ds),\\
 \int_0^t\int_Z \biggl(\int_0^1 \Phi^{\prime \prime}(Y(s-)+\theta \Psi)\Psi^2 d\theta \biggr)\tilde{\eta}(dz,ds),
\end{align*}
are martingales with zero mean. Hence taking the mathematical expectation yields
\begin{equation}\label{inv3}
\begin{split}
\EE \Phi(Y(t)-\Phi(\lvert\xi\rvert^2))
 \le  \EE \int_0^t \Phi^\prime (Y(s)) \biggl(\int_Z \lvert \sigma(s, \bu(s;\xi), z) \rvert^2 \nu(dz)\biggr)ds\\
  +\EE \int_0^t\int_Z \biggl(\int_0^1 \Phi^{\prime \prime}(Y(s-)+\theta \Psi)\Psi^2 d\theta \biggr)\nu(dz)ds\\
  -2\kappa_1 \EE \int_0^t \Phi^\prime (Y(s))\lVert \bu(s;\xi)\rVert^2_2 ds.
\end{split}
\end{equation}
Since
\begin{align*}
 \Phi^{\prime\prime}(Y(s-)+\theta \Psi)\Psi^2=&\frac{-2\eps \Psi^2}{(1+\eps Y(s)+\eps \theta \Psi)^3}\\
 =&\frac{-2\eps \Psi^2}{(1+\eps \theta\lvert \sigma(s, \bu(s,\xi), z)+\bu(s-;\xi)\rvert^2-\eps(1-\theta)\lvert\bu(s-;\xi)\rvert^2)},
\end{align*}
we see that $\Phi^{\prime\prime}(Y(s-)+\theta \Psi)\Psi^2 \le 0$ for any $\theta \in [0,1]$. Therefore we can drop out the second term in the right-hand side of
\eqref{inv3}, use item (2) in Condition \ref{cond-s-1} to obtain
\begin{equation}\label{inv4-b}
\begin{split}
 \EE \Phi(Y(t))+2\kappa_1 \EE \int_0^t \Phi^\prime (Y(s))\lVert \bu(s;\xi)\rVert^2_2 ds
  \le \Phi(\lvert\xi\rvert^2)  + \ell_1 \EE \int_0^t \Phi^\prime (Y(s)) \lvert \bu(s; \xi)\rvert^2ds\\
  + \ell_0 \EE \int_0^t \Phi^\prime (Y(s)) ds.
\end{split}
\end{equation}
By using Poincar\'e's inequality (see \eqref{POINCARE}) the last estimate becomes
\begin{equation}\label{inv5}
\begin{split}
 \EE \Phi(Y(t))+2\kappa_1 \lambda_1^2 \EE \int_0^t \Phi^\prime (Y(s))\lvert \bu(s;\xi)\rvert^2_2 ds
  \le \Phi(\lvert\xi\rvert^2)  + \ell_1 \EE \int_0^t \Phi^\prime (Y(s)) \lvert \bu(s; \xi)\rvert^2ds\\
  + \ell_0 \EE \int_0^t \Phi^\prime (Y(s)) ds.
\end{split}
\end{equation}
By integrating both side of this last inequality wrt $\mu$ on $H$ and using the fact that
\begin{equation}\label{inv5-b}
\int_H \EE \phi(\bu(s,\xi))\mu(dx)=\int_H \phi(\xi) \mu(dx), \forall \phi\in C_b(H), \text{  ($\mu$ is an invariant measure)}
\end{equation}

we obtain from \eqref{inv5} that
\begin{equation*}
 \left(2\kappa_1\lambda_1^2-\ell_1\right)\int_H \frac{\lvert \xi\rvert^2}{(1+\eps \lvert \xi\rvert^2)^2}\mu(dx)\le \ell_0 \int_H \frac{1}{(1+\eps \lvert \xi\rvert^2)^2}\mu(dx).
\end{equation*}
From this inequality we obtain that
\begin{equation}\label{inv6}
 \int_H \frac{\lvert \xi\rvert^2}{(1+\eps \lvert \xi\rvert^2)^2}\mu(dx)\le \frac{\ell_0}{2\kappa_1 \lambda_1^2-\ell_1},
\end{equation}
where we have used the facts that $\frac{1}{(1+\eps \lvert \xi\rvert^2)^2}\ge 1$, $2\kappa_1\lambda_1^2-\ell_1>0$  and $\mu(V)+\mu(H\backslash V)=1$.

From \eqref{inv4-b} and \eqref{inv6} we derive that
\begin{equation}
 2\kappa_1\int_H \EE\int_0^t \Phi^\prime(Y(s))\lVert \bu(s;\xi)\rVert^2_2 ds\le \frac{\ell_0}{2\kappa_1 \lambda_1^2-\ell_1}\left(\ell_1+1\right)+\ell_0.
\end{equation}
Choosing $\phi(\bu(s;\xi))=\int_0^t \Phi^\prime (Y(s))\lVert\bu(s;\xi)\rVert^2_2 ds$ and using \eqref{inv5-b} we see that
\begin{equation}\label{inv7}
 \int_H \frac{\lVert \xi\rVert^2_2}{(1+\eps \lvert \xi\rvert^2)^2} \mu(dx)\le \frac{\ell_0}{2\kappa_1(2\kappa_1\lambda_1^2-\ell_1)}\left(\ell_1+1\right)+\frac{\ell_0}{2\kappa_1}.
\end{equation}
Adding up \eqref{inv6} and \eqref{inv7} side by side, letting $\eps\rightarrow 0$ and using Fatou's lemma imply that
\begin{equation}\label{inv8}
 \int_H \left(\lvert \xi\rvert^2+\lVert \xi\rVert^2_2\right)\mu(dx)\le \frac{\ell_0}{2\kappa_1\lambda_1^2-\ell_1}\biggl(\frac{\ell_1+1}{2\kappa_1}+1\biggr)+\frac{\ell_0}{2\kappa_1},
\end{equation}
which terminates the proof of the proposition.
\end{proof}
We can prove the ergodicity of the invariant measure under the condition that $\kappa_1$ is large enough.
\begin{thm}
Assume that $2\kappa_1\lambda_1^2>\ell_1$. Then, the Markovian semigroup $\CP_t$ has an invariant measure $\mu$ which is tight  and ergodic on $H$.
\end{thm}
\begin{proof}
 Let $\mathcal{M}\subset \mathcal{M}_1(H)$ be the set of invariant measure of $\CP_t$ and
 $$\tilde{\ell}=\frac{\ell_0}{2\kappa_1\lambda_1^2-\ell_1}\biggl(\frac{\ell_0}{2\kappa_1}+1\biggr)+\frac{\ell_0}{2\kappa_1}.$$ It is not difficult to show that
 $\mathcal{M}$ is convex (see for example \cite[page 296]{Hasselblatt}). As before let $R>0$ and $A_R=\{\bu\in H: \lVert \bu\rVert_2>R\}$. We see from Chebychev-Markov's inequality that
 \begin{equation*}
  \sup_{\mu\in \mathcal{M}}\mu\left(A_R\right)\le \frac{1}{R^2}\int_H\lVert \xi\rVert^2_2 \mu(dx).
 \end{equation*}
Owing to \eqref{inv7} we have that
\begin{equation*}
  \sup_{\mu\in \mathcal{M}}\mu\left(A_R\right)\le \frac{\tilde{\ell}}{R^2},
\end{equation*}
which implies that for any $\eps>0$
\begin{equation*}
 \mu(B_{V}(\frac{1}{\sqrt{\eps}}))\ge 1-\eps,
\end{equation*}
where $B_{V}(\frac{1}{\sqrt{\eps}})=V\backslash A_{\frac{1}{\sqrt{\eps}}}$. Since $B_{V}(\frac{1}{\sqrt{\eps}})$
is compact in $H$ we infer that the set $\mathcal{M}$ is tight on $H$.
Since $\mathcal{M}$ is non-empty, convex and tight, by Krein-Millman's theorem (see, for instance, \cite[Theorem 3.65, p. 110]{Fabian}) it has extrema which are ergodic.
We deduce from the above argument that $\CP_t$ has at least one invariant measure which is ergodic.
\end{proof}

\section*{Acknowledgment}
The authors' research is supported by the Austrian Science Foundation through the
grant number P20705.


\begin{thebibliography}{99}

\bibitem{Alb+Brz+Wu_2010} S. Albeverio, Z.  Brze\'zniak, J.L. Wu,
\textit{Existence of global solutions and invariant measures for stochastic differential equations driven by Poisson type noise
with non-Lipschitz coefficients},
J. Math. Anal. Appl. \textbf{371}, no. 1, 309-322  (2010).

% \bibitem{Alb+Brz+Wu_2010} S.~Albeverio, Z.~Brzezniak and J.-L.~Wu. Existence of
% global solutions and invariant measures for stochastic differential
% equations driven by Poisson type noise with non-Lipschitz coefficients.
% \newblock{\em J. Math. Anal. Appl. } 371(1):309--322, 2010.

\bibitem{Albeverio+Mandrekar+Rudiger}S.~Albeverio, V.~Mandrekar and B.~ R\"udiger.\newblock{
Existence of mild solutions for stochastic differential equations
and semilinear equations with non-Gaussian L\'evy noise.}
\newblock{\em Stochastic Process. Appl.} 119(3):835-863, 2009.

\bibitem{BABIN+VISHIK}A. V.~Babin and M. I.~Vishik. \newblock{ \em Attractors of evolution equations.}
Volume 25 of \newblock{\em
 Studies in Mathematics and its Applications}. North-Holland Publishing Co., Amsterdam, 1992.





\bibitem{BELLOUT2} H.~Bellout, F.~Bloom and J.~Necas. %
\newblock{Phenomenological behavior of multipolar viscous fluids.} %
\newblock{\em Quarterly of Applied Mathematics} 50:559-583, 1992.

\bibitem{BELLOUT1} H.~Bellout, F.~Bloom and J.~Necas. \newblock{Solutions
for incompressible Non-Newtonian fluids}. \newblock{\em C. R. Acad. Sci.
Paris S\'er I. Math.} 317:795-800,1993.

\bibitem{BELLOUT4} H.~Bellout, F.~Bloom and J.~Necas. \newblock{Young
measure-valued solutions for Non-Newtonian incompressible fluids.} %
\newblock{Communication in Partial Differential Equations.} 19(11\&
12):1763-1803, 1994.

\bibitem{BELLOUT3} H.~Bellout, F.~Bloom and J.~Necas. \newblock{Bounds for
the dimensions of the attractors of nonlinear bipolar viscous fluids}. %
\newblock{\em Asymptotic Analysis.} 11(2):131-167,1995.

\bibitem{BELLOUT5} H.~Bellout, F.~Bloom and J.~Necas. \newblock{Existence,
uniqueness and stability of solutions to initial boundary value problems for
bipolar fluids.}  \newblock{\em Differential and Integral Equations}
8:453-464, 1995



\bibitem{bensoussan} A.~Bensoussan. \newblock Stochastic {N}avier-{S}tokes {E%
}quations. \newblock {\em Acta Applicandae Mathematicae}, 38:267--304, 1995.

\bibitem{Bensoussan-Temam} A.~Bensoussan and R.~Temam. \newblock Equations {S%
}tochastiques du {T}ype {N}avier-{S}tokes.
\newblock {\em Journal of
Functional Analysis}, 13:195--222, 1973.


\bibitem{breckner} H.~Breckner
\newblock{\em Approximation and optimal
control of the stochastic navier-Stokes equation}
\newblock{Dissertation,
Martin-Luther University, Halle-Wittenberg}, 1999.

\bibitem{Breckner} H.~Breckner \newblock{\em Galerkin approximation and the
strong solution of the Navier-Stokes equation}. \newblock{\em Journal of
Applied Mathematics and Stochastic Analysis}. 13(3):239--259, 2000.

\bibitem{Brzezniak3} Z.~Brzezniak and L.~Debbi. On stochastic Burgers
equation driven by a fractional Laplacian and space-time white noise. %
\newblock{\em Stochastic differential equations: Theory and applications,}
Interdiscip. Math. Sci., 2, World Sci. Publ., pages 135--167, 2007.

\bibitem{Brz+Haus_2009} Z. Brze\'{z}niak, E. Hausenblas. \newblock{Maximal regularity for stochastic convolutions driven by L\'evy processes},
 \newblock{\em Probab. Theory Related Fields.}  \textbf{145}(3-4):615--637, 2009.

\bibitem{ZB+EH+JZ}Z.~Brze\'{z}niak, E.~Hausenblas and J.~Zhu.{2D stochastic
Navier-Stokes equations driven by jump noise.} \newblock{Preprint.}

\bibitem{Brzezniak4} Z.~Brzezniak, B.~Maslowski and J.~Seidler. Stochastic
nonlinear beam equation. \newblock{\em Probab. Theory Relat. Fields}.
132(2):119--144, 2005.

\bibitem{Caraballo1} T.~Caraballo, J.A.~Langa and T.~Taniguchi. The
exponential behaviour and stabilizability of stochastic 2D-Navier-Stokes
equations. \newblock{\em J. Differential Equations.} 179(2):714--737, 2002.

\bibitem{Caraballo} T.~Caraballo, J.~Real and T.~Taniguchi. {{On the
existence and uniqueness}} of solutions to stochastic three-dimensional
Lagrangian averaged Navier-Stokes equations. \newblock{\em Proc. R. Soc.
Lond. Ser. A Math. Phys. Eng. Sci.} 462(2066):459--479, 2006.

\bibitem{Caraballo3} T.~Caraballo, A.M.~{{M\'arquez-Dur\'an}} and
J.~Real. \newblock{The asymptotic behaviour of a stochastic 3D LANS-$\alpha$
model.} \newblock{\emph Appl. Math. Optim.} 53(2):141--161, 2006.

% \bibitem{CHEN} J.~Chen and Z.M.~Chen. \newblock{Stochastic Non-Newtonian
% fluid motion equations of a nonlinear bipolar viscous fluid.} \newblock{\em
% Journal of Mathematical Analysis and Applications.} 369:486-509, 2010.

\bibitem{Chow+Kashminskii}P.-L.~Chow and R.~Z.~Khasminskii.\newblock{
Stationary solutions of nonlinear stochastic evolution equations.}
\newblock{\em Stochastic Anal. Appl.} 15(5):671-699, 1997.

\bibitem{Constantin+Foias+Temam} P.~Constantin, C.~Foias and R.~Temam.
Attractors representing turbulent flows.
\newblock{\em Mem. Amer. Math. Soc.} Vol 53, no. 314, vii+67 pp, 1985.

\bibitem{Constantin+Foias+Temam-2}P.~Constantin, C.~Foias and R.~Temam. \newblock{.}
\newblock{\em Physica D: Nonlinear Phenomena} 30(3):284-296, 1988.

\bibitem{Millet+Chueshov_2010} I. Chueshov and A.  Millet,
 \newblock{Stochastic 2D hydrodynamical type systems: well posedness and large deviations}, \newblock{\em Appl. Math. Optim.} \textbf{61}(3):379-420, 2010.


\bibitem{daprato} G.~Da~Prato and J.~Zabczyk.
\newblock {\em Stochastic
{E}quations in {I}nfinite {D}imensions}. \newblock Cambridge University
Press, 1992.

\bibitem{DAPRATO} G.~Da Prato and A.~Debussche {2D stochastic Navier-Stokes
equations with a time-periodic forcing term.} \emph{J. Dynam. Differential
Equations} 20(2):301--335, 2008.

% \bibitem{DEUGOUE1} G.~Deugoue and M.~Sango. On the Stochastic 3D
% Navier-Stokes-$\alpha$ Model of Fluids Turbulence,  \textit{Abstract and
% Applied Analysis}, vol. 2009, Article ID 723236, 27 pages, 2009.
% doi:10.1155/2009/723236.

\bibitem{DEUGOUE2} G.~Deugoue and M.~Sango. On the Strong Solution for the
3D Stochastic Leray-Alpha Model, \textit{Boundary Value Problems}, vol.
2010, Article ID 723018, 31 pages, 2010. doi:10.1155/2010/723018.

\bibitem{DongNS} Z. Dong and Z.~Jianliang. {Martingale solutions and Markov selection of stochastic 3D
Navier-Stokes equations with jump}, {\em  J. Differential
Equations}, \textbf{250}:27372778, 2011.

\bibitem{Hasselblatt} B.~Hasselblatt and A.~Katok. {\em A first course in dynamics. With a panorama of recent developments.} Cambridge University Press, New York, 2003.

\bibitem{EH+PAR}E.~Hausenblas, P.~A.~Razafimandimby and M.~ Sango.
Martingale solution to Differential type fluids of grade two driven by random force of L\'evy type. \newblock{\em Submitted.}

\bibitem{Fabian}M.~Fabian, P.~Habala, P.~Hájek, V.~Montesinos and V.~Zizler. {\em Banach space theory. The basis for linear and nonlinear analysis.}
CMS Books in Mathematics/Ouvrages de Math\'ematiques de la SMC. Springer, New York, 2011.

% \bibitem{DU} Q.~Du and M.D.~Gunzburger. \newblock{Analysis of Ladyzhenskaya
% Model for Incompressible Viscous Flow.}  \newblock{\em Journal of
% Mathematical Analysis and Applications.} 155:21-45, 1991.

\bibitem{Flandoli2} F.~Flandoli, M.~Gubinelli, M.~Hairer, and M.~Romito.
\newblock{Rigourous remarks about scaling laws in turbulent
fluid.} \newblock{\em Commun. Math. Phys.} 278: 1-29, 2008.

\bibitem{RUZICKA} J.~Freshe and M.~Ruzicka. \newblock{Non-homogeneous
generalized Newtonian fluids.} \newblock{Mathematische Zeitschrift}.
260(2):353-375, 2008.

\bibitem{Gyongy-Krylov} I. Gy\"{o}ngy and N.V. Krylov. On stochastics equations with respect to semimartingales.
 II. It\^o formula in Banach spaces. \newblock{\em Stochastics.} 6(3-4):153--173, 1981/1982.

% \bibitem{GUO2} B.~Guo and G.~Lin. \newblock{Existence and Uniqueness of
% stationary solutions of Non-Newtonian viscous incompressible fluids.} %
% \newblock{\em  Communications in nonlinear Science and Numerical Simulation.}
% 4(1):63-68, 1999.
%
\bibitem{GUO1} B.~Guo, C.~Guo and J.~Zhang. \newblock{Martingale and
stationary solutions for stochastic Non-Newtonian fluids.}  \newblock{\em
Differential and Integral Equations.} 23(3 \& 4):303-326, 2010.

\bibitem{Ik-Wat-81}
N.~Ikeda and S.~Watanabe,
\newblock {em Stochastic differential equations and diffusion processes},
  volume~ \textbf{24} of {North-Holland Mathematical Library}.
\newblock North-Holland Publishing Co., Amsterdam, second edition, 1989.

\bibitem{KAP} A.~Kupiainen. Statistical theories of turbulence.
In \newblock{\em Advances in Mathematical Sciences and Applications.} Gakkotosho, Tokyo, 2003.



\bibitem{LADY1} O.A.~Ladyzhensakya. \newblock{\em The mathematical theory of
viscous incompressible flow.} Gordon and Breach, New York, 1969.

\bibitem{LADY2} O.A.~Ladyzhensakya. \newblock{New equations for the
description of the viscous incompressible fluids and solvability in the
large of the boundary value problems for them.}  In \newblock{\em Boundary
Value Problems of Mathematical Physics V.} American Mathematical Society,
Providence, RI, 1970.
%
% \bibitem{LIONS} J.L.~Lions. \newblock{\em Quelques m\'ethodes de
% r\'esolution de probl\'ems aux limites non-lin\'eaires.} Dunod, Paris, 1969.

\bibitem{MALEK} J.~Malek, J.~Necas and A.~Novotny. \newblock{Measure-valued
solutions and asymptotic behavior of a multipolar model of a boundary layer.}
\newblock{\em Czechoslovak Mathematical Journal.} 42(3):549-576, 1992.

\bibitem{ROKYTA} J.~Malek, J.~Necas, M.~Rokyta and M.~Ruzicka.
\newblock{\em Weak and measure-valued solutions to evolutionary
PDEs.} Applied Mathematics and Mathematical Computation, 13. Chapman \&
Hall, London, 1996.

\bibitem{ROZOVSKII1} R.~Mikulevicius and B.L.~ Rozovskii. %
\newblock{Stochastic {N}avier-{S}tokes {E}quations and {T}urbulent {F}lows.} %
\newblock{\em SIAM J. Math. Anal.}, 35(5):1250-1310, 2004.



\bibitem{NECAS2} J.~Necas and M.~Silhavy \newblock{Multipolar viscous
fluids.} \newblock{\em Quaterly of Applied Mathematics.} XLIX(2):247-266,
1991.

\bibitem{NECAS1} J.~Necas, A.~Novotny and M.~Silhavy. \newblock{Global
solution to the compressible isothermal multipolar fluids}. \newblock{\em J.
Math. Anal. Appl.} 162:223-242, 1991.


\bibitem{pardoux} E.~Pardoux. \emph{Equations aux d\'eriv\'ees partielles
stochastiques monotones}. Th\`ese de {D}octorat, Universit\'e Paris-Sud,
1975.

\bibitem{Peszat} S.~Peszat and J.~Zabczyk.
\newblock{\em Stochastic Partial Differential Equations with Levy Noise. An evolution
equation approach}. Encyclopedia of Mathematics and its Applications 113,
Cambridge university Press, 2007.

\bibitem{PAUL1} P. A. Razafimandimby and M. Sango. Weak Solutions of a
Stochastic Model for Two-Dimensional Second Grade Fluids,\textit{Boundary
Value Problems}. vol. 2010, Article ID 636140, 47 pages, 2010.
doi:10.1155/2010/636140.

\bibitem{PAUL2} P. A. Razafimandimby. On Stochastic Models Describing the
Motions of Randomly Forced Linear Viscoelastic Fluids, \textit{Journal of
Inequalities and Applications}, vol. 2010, Article ID 932053, 27 pages,
2010. doi:10.1155/2010/932053.


\bibitem{PAUL3} P. A.~Razafimandimby and M.~Sango.
\newblock{Asymptotic
behavior of solutions of stochastic evolution equations for second
grade fluids.} \newblock{\em C. R. Math. Acad. Sci. Paris}, Volume 348,
Issues 13-14, Pages 787-790, 2010.

\bibitem{PAUL-13} P.~A.~Razafimandimby and M.~Sango.
  Strong solution for a stochastic model of two-dimensional second grade fluids: existence, uniqueness and asymptotic behaviour.
  \newblock{\em Nonlinear Analysis-Theory Methods \& Applications.}  75(11):4251--4270, 2012.

\bibitem{Robinson}J.~C.~Robinson. \newblock{Infinite-Dimensional Dynamical Systems: An Introduction to Dissipative Parabolic PDEs and the Theory of Global Attractors,}
Cambridge Texts in Applied Mathematics, Cambridge University Press, Cambridge, 2001.

\bibitem{Rudiger} B.~R\"udiger. Stochastic integration with respect to compensated Poisson random measure on separable Banach spaces.
\newblock{Stochastic and Stochastic Reporst.} 76(3): 213-242, 2004.
%
% \bibitem{sango1} M.~Sango. Weak solutions for a doubly degenerate
% quasilinear parabolic equation with random forcing. \emph{Discrete Contin.
% Dyn. Syst. Ser. B}, 7(4) (2007) 885--905.

\bibitem{sango2} M.~Sango. Magnetohydrodynamic turbulent flows: Existence
results.\textit{\ Physica D: Nonlinear Phenomena} 239(12): 912-923, 2010.

\bibitem{sango} M.~Sango. Density dependent stochastic Navier-Stokes
equations with non Lipschitz random forcing.  \textit{Reviews in
Mathematical Physics} 22(6):669--697, 2010.

\bibitem{Taniguchi}T.~Taniguchi. \newblock{The existence and asymptotic behaviour of energy solutions to stochastic
2D functional Navier-Stokes equations driven by Levy processes.}
\newblock{\em Journal of Mathematical Analysis and Applications}. 385(2):634-654, 2012.



\bibitem{temam} R.~Temam. \newblock {\em Navier-{S}tokes {E}quations}. %
\newblock North-Holland, 1979.

\bibitem{Temam-Infinite}R.~Temam. \newblock{\em Infinite-dimensional dynamical systems in mechanics and physics.} In: Applied Mathematical Sciences,
vol. 68. Springer-Verlag, New York, 1988.

\bibitem{Zhu_2010} J. Zhu, {\em A Study of SPDEs w.r.t. Compensated Poisson Random Measures and Related Topics}, Ph. D. Thesis,  University of York, 2010.
\end{thebibliography}
\end{document}